\documentclass[aps,12pt,prd,citeautoscript,longbibliography,superscriptaddress,nofootinbib]{revtex4-1}
\usepackage[T1]{fontenc}
\usepackage{graphicx}
\usepackage{amsmath,amssymb,amsfonts,bm,dsfont,units}
\usepackage{mathtools}
\usepackage{hyperref}

\usepackage{yfonts}

\newcommand{\coho}[1]{\textswab{#1}}

\usepackage[mathscr]{euscript}

\usepackage{braket}

\usepackage{color}
\definecolor{dkgreen}{rgb}{0,0.5,0}
\definecolor{midnightblue}{rgb}{0.39,0.58,0.93}
\definecolor{ltgreen}{rgb}{0.1,.59,.43}
\definecolor{hanpurple}{rgb}{0.32, 0.09, 0.98}
\definecolor{readableyellow}{rgb}{0.55, 0.45,0}
\hypersetup{colorlinks=true,citecolor=ltgreen,linkcolor=ltgreen,urlcolor=hanpurple,pdfstartview=FitH}

\DeclareMathAlphabet{\mathpzc}{OT1}{pzc}{m}{it}

\definecolor{purpledk}{rgb}{0.32, 0.09, 0.80}

\newcommand{\I}{\mathcal{I}}

\newcommand{\spt}{\operatorname{\bf SPT}}

\newcommand{\bMTC}{\mathcal{C}}

\newcommand{\defectO}{\mathscr{O}}

\newcommand{\Xdot}{\mathscr{X}}

\newcommand{\rel}{\mathscr{R}}

\newcommand{\cocyleatorR}{\mathord{\vcenter{\hbox{\includegraphics[scale=1]{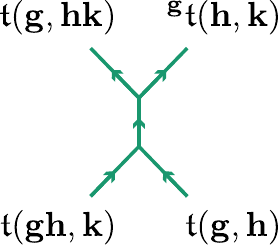}}}}}
\newcommand{\cocyleatorL}{\mathord{\vcenter{\hbox{\includegraphics[scale=1]{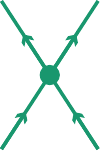}}}}}

\newcommand{\FMoveCocyR}{\mathord{\vcenter{\hbox{\includegraphics[scale=1]{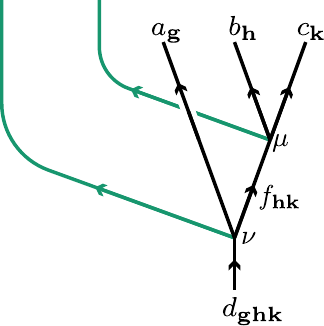}}}}}
\newcommand{\FMoveCocyL}{\mathord{\vcenter{\hbox{\includegraphics[scale=1]{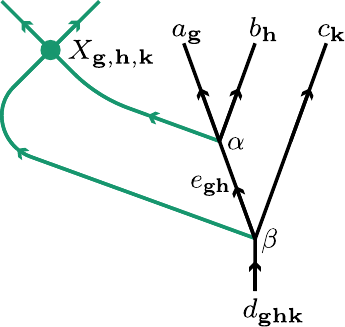}}}}}
\newcommand{\FMoveCocyleLBranched}{\mathord{\vcenter{\hbox{\includegraphics[scale=1]{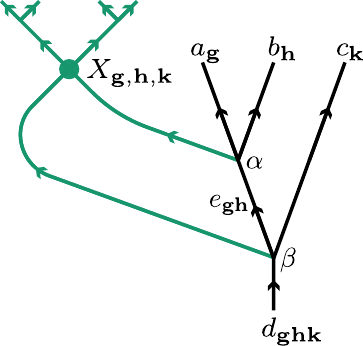}}}}}

\newcommand{\RelabelLeftF}{\mathord{\vcenter{\hbox{\includegraphics[scale=1]{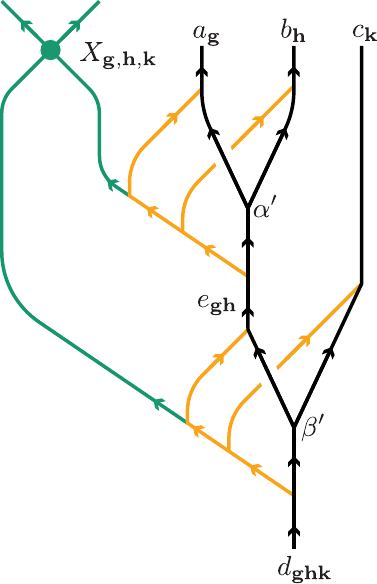}}}}}
\newcommand{\RelabelRightF}{\mathord{\vcenter{\hbox{\includegraphics[scale=1]{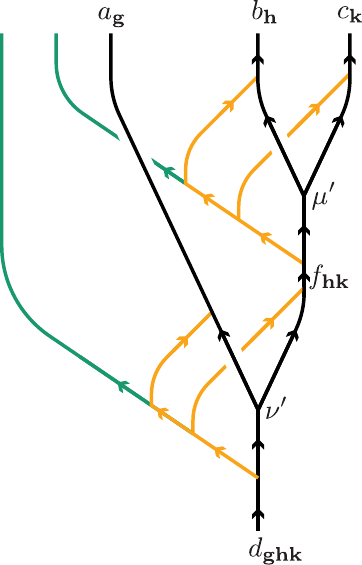}}}}}

\newcommand{\XConsisL}{\mathord{\vcenter{\hbox{\includegraphics[scale=1]{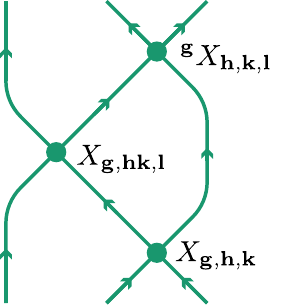}}}}}
\newcommand{\XConsisR}{\mathord{\vcenter{\hbox{\includegraphics[scale=1]{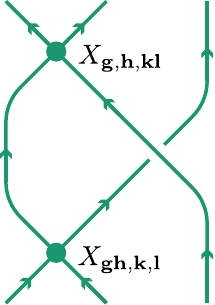}}}}}

\newcommand{\FusionSpaceL}{\mathord{\vcenter{\hbox{\includegraphics[scale=1]{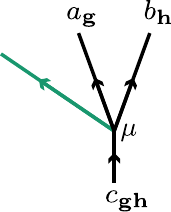}}}}}
\newcommand{\FusionSpaceR}{\mathord{\vcenter{\hbox{\includegraphics[scale=1]{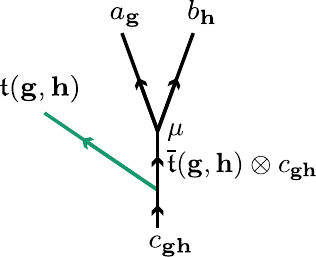}}}}}

\newcommand{\FusionSpaceHat}{\mathord{\vcenter{\hbox{\includegraphics[scale=1]{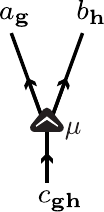}}}}}

\newcommand{\FusionSpaceLDual}{\mathord{\vcenter{\hbox{\includegraphics[scale=1]{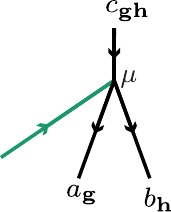}}}}}
\newcommand{\FusionSpaceRDual}{\mathord{\vcenter{\hbox{\includegraphics[scale=1]{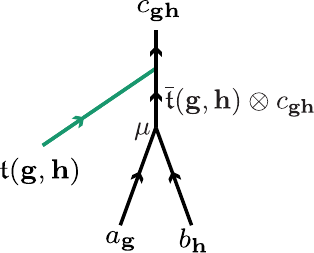}}}}}
\newcommand{\FusionSpaceHatdual}{\mathord{\vcenter{\hbox{\includegraphics[scale=1]{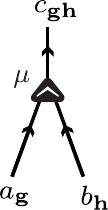}}}}}

\newcommand{\RHatResolve}{\mathord{\vcenter{\hbox{\includegraphics[scale=1]{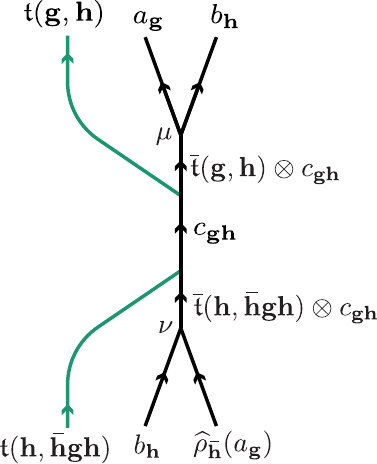}}}}}

\newcommand{\RHatL}{\mathord{\vcenter{\hbox{\includegraphics[scale=1]{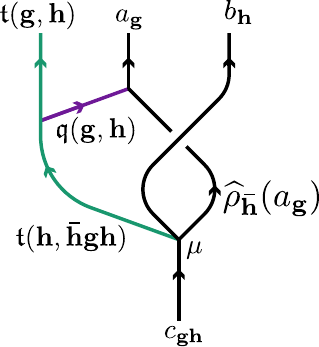}}}}}
\newcommand{\RHatR}{\mathord{\vcenter{\hbox{\includegraphics[scale=1]{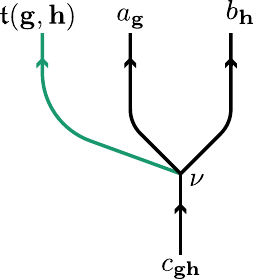}}}}}

\newcommand{\RRelabelL}{\mathord{\vcenter{\hbox{\includegraphics[scale=1]{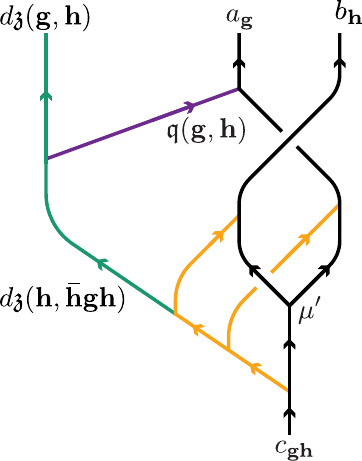}}}}}
\newcommand{\RRelabelR}{\mathord{\vcenter{\hbox{\includegraphics[scale=1]{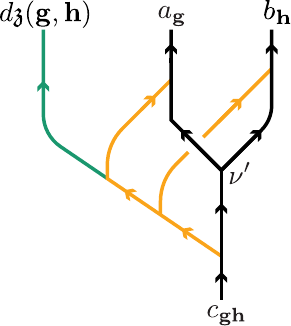}}}}}

\newcommand{\RHatOp}{\mathord{\vcenter{\hbox{\includegraphics[scale=1]{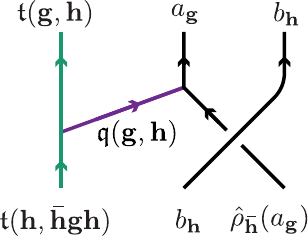}}}}}
\newcommand{\RHatOpInv}{\mathord{\vcenter{\hbox{\includegraphics[scale=1]{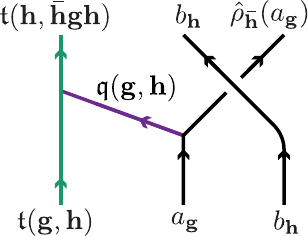}}}}}

\newcommand{\QLeft}{\mathord{\vcenter{\hbox{\includegraphics[scale=1]{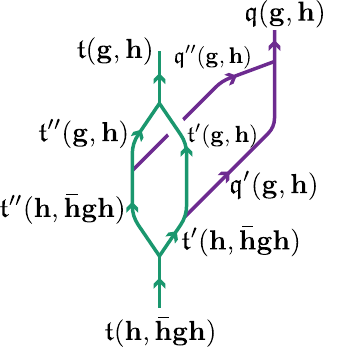}}}}}
\newcommand{\QRight}{\mathord{\vcenter{\hbox{\includegraphics[scale=1]{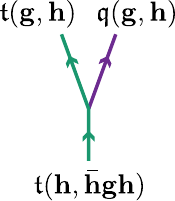}}}}}

\newcommand{\PUdef}{\mathord{\vcenter{\hbox{\includegraphics[scale=1]{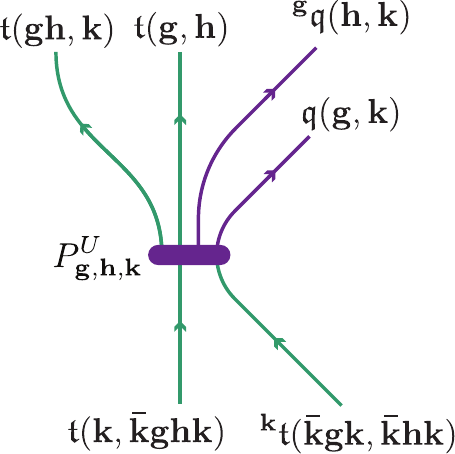}}}}}
\newcommand{\PUdefRHS}{\mathord{\vcenter{\hbox{\includegraphics[scale=1]{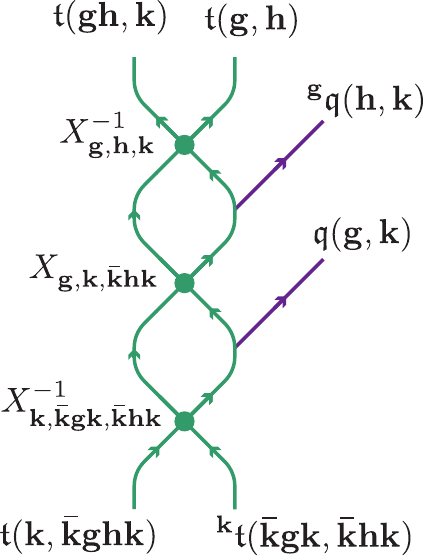}}}}}

\newcommand{\PEtaDef}{\mathord{\vcenter{\hbox{\includegraphics[scale=1]{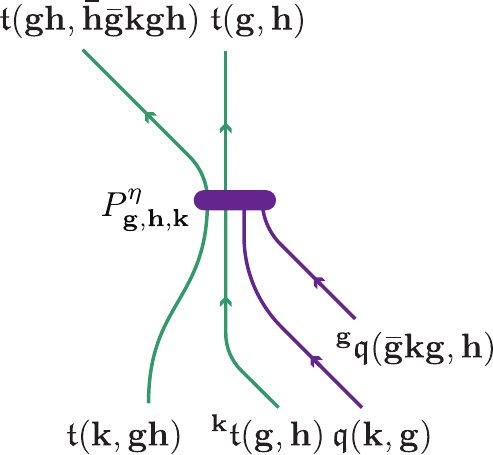}}}}}
\newcommand{\PEtaDefR}{\mathord{\vcenter{\hbox{\includegraphics[scale=1]{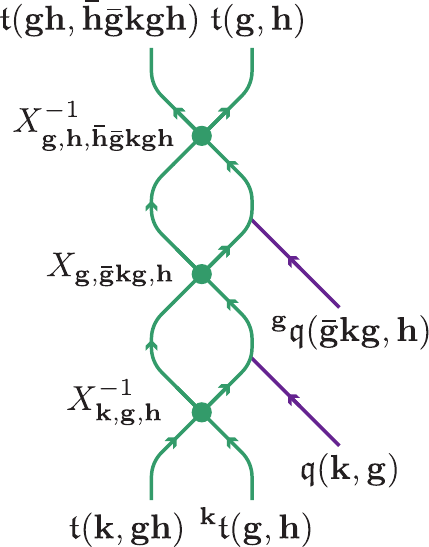}}}}}

\newcommand{\FLeft}{\mathord{\vcenter{\hbox{\includegraphics[scale=1]{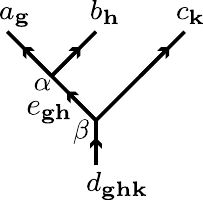}}}}}
\newcommand{\FRight}{\mathord{\vcenter{\hbox{\includegraphics[scale=1]{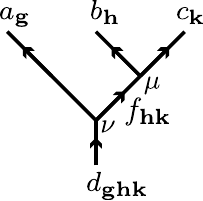}}}}}

\newcommand{\RLeft}{\mathord{\vcenter{\hbox{\includegraphics[scale=1]{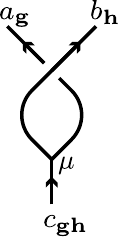}}}}}
\newcommand{\RRight}{\mathord{\vcenter{\hbox{\includegraphics[scale=1]{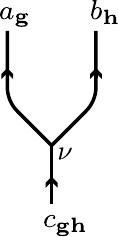}}}}}

\newcommand{\RLeftHat}{\mathord{\vcenter{\hbox{\includegraphics[scale=1]{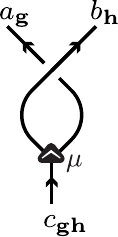}}}}}
\newcommand{\RRightHat}{\mathord{\vcenter{\hbox{\includegraphics[scale=1]{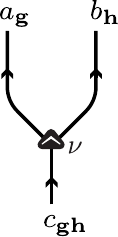}}}}}

\newcommand{\ULeft}{\mathord{\vcenter{\hbox{\includegraphics[scale=1]{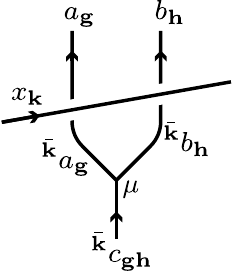}}}}}
\newcommand{\URight}{\mathord{\vcenter{\hbox{\includegraphics[scale=1]{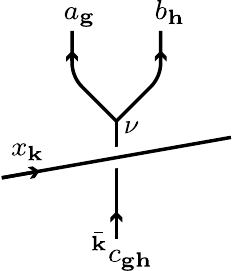}}}}}

\newcommand{\etaLeft}{\mathord{\vcenter{\hbox{\includegraphics[scale=1]{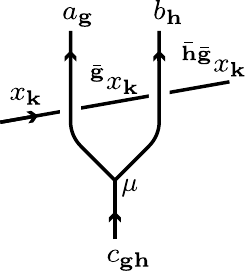}}}}}
\newcommand{\etaRight}{\mathord{\vcenter{\hbox{\includegraphics[scale=1]{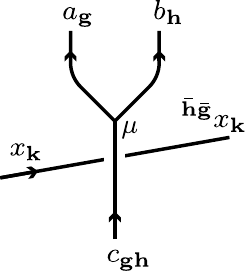}}}}}

\newcommand{\UHatLPrime}{\mathord{\vcenter{\hbox{\includegraphics[scale=1]{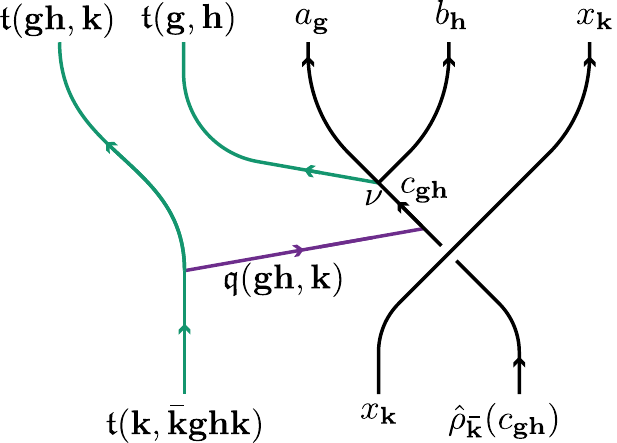}}}}}
\newcommand{\UHatRPrime}{\mathord{\vcenter{\hbox{\includegraphics[scale=1]{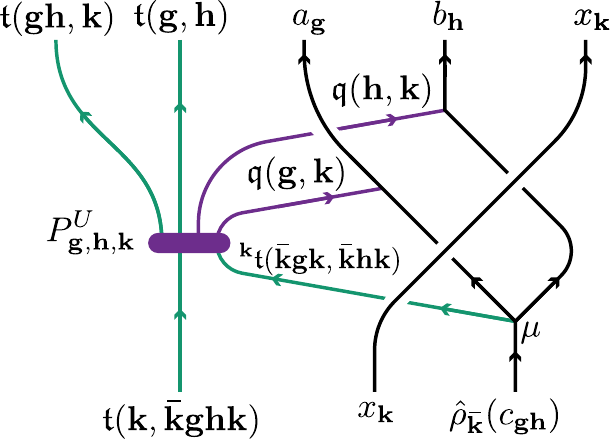}}}}}

\newcommand{\etaHatLPrime}{\mathord{\vcenter{\hbox{\includegraphics[scale=1]{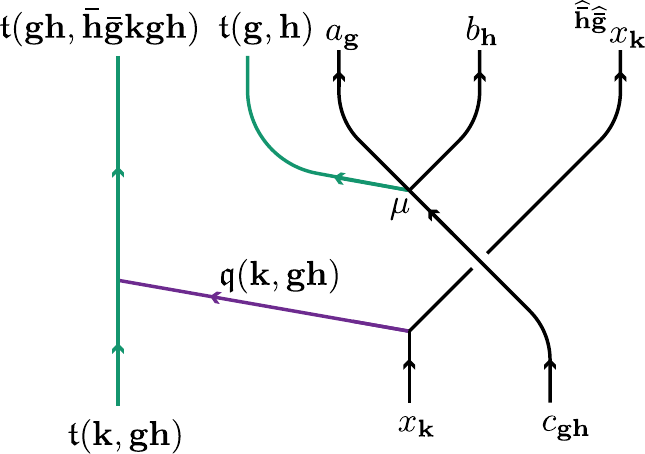}}}}}
\newcommand{\etaHatRPrime}{\mathord{\vcenter{\hbox{\includegraphics[scale=1]{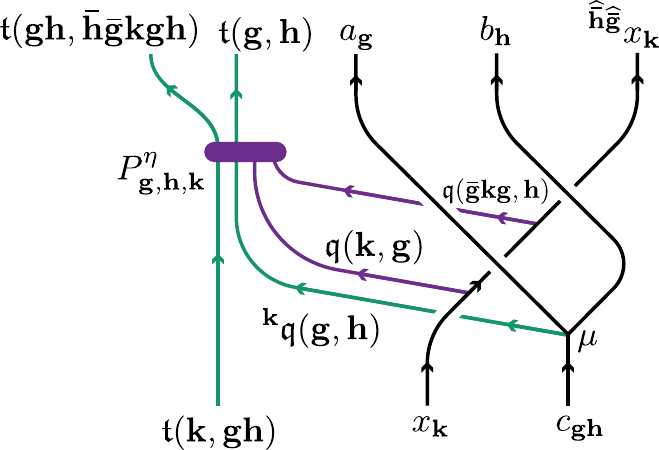}}}}}

\newcommand{\VabcmuL}{\mathord{\vcenter{\hbox{\includegraphics[scale=1]{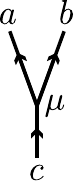}}}}}
\newcommand{\VabcmuR}{\mathord{\vcenter{\hbox{\includegraphics[scale=1]{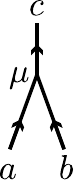}}}}}

\newcommand{\RVabc}{\mathord{\vcenter{\hbox{\includegraphics[scale=1]{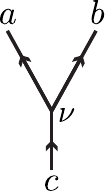}}}}}
\newcommand{\Rabc}{\mathord{\vcenter{\hbox{\includegraphics[scale=1]{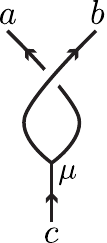}}}}}

\newcommand{\Rab}{\mathord{\vcenter{\hbox{\includegraphics[scale=1]{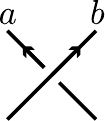}}}}}
\newcommand{\Mab}{\mathord{\vcenter{\hbox{\includegraphics[scale=1]{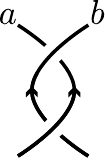}}}}}
\newcommand{\idab}{\mathord{\vcenter{\hbox{\includegraphics[scale=1]{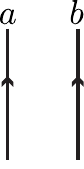}}}}}

\newcommand{\idabesolve}{\mathord{\vcenter{\hbox{\includegraphics[scale=1]{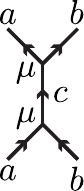}}}}}
\newcommand{\idc}{\mathord{\vcenter{\hbox{\includegraphics[scale=1]{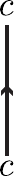}}}}}
\newcommand{\bubbleabc}{\mathord{\vcenter{\hbox{\includegraphics[scale=1]{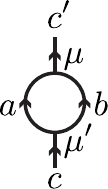}}}}}

\newcommand{\FLeftMTC}{\mathord{\vcenter{\hbox{\includegraphics[scale=1]{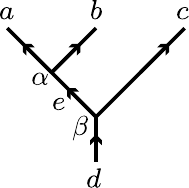}}}}}
\newcommand{\FRightMTC}{\mathord{\vcenter{\hbox{\includegraphics[scale=1]{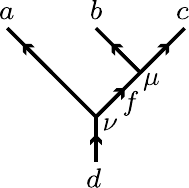}}}}}

\newcommand{\FLeftHat}{\mathord{\vcenter{\hbox{\includegraphics[scale=1]{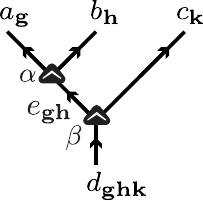}}}}}
\newcommand{\FRightHat}{\mathord{\vcenter{\hbox{\includegraphics[scale=1]{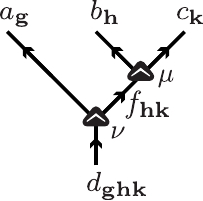}}}}}

\newcommand{\braidingop}{\mathord{\vcenter{\hbox{\includegraphics[scale=1]{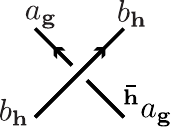}}}}}

\newcommand{\relabelbosonicR}{\mathord{\vcenter{\hbox{\includegraphics[scale=1]{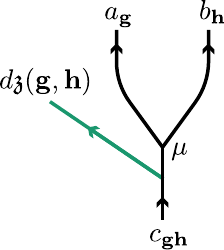}}}}}
\newcommand{\relabelbosonicL}{\mathord{\vcenter{\hbox{\includegraphics[scale=1]{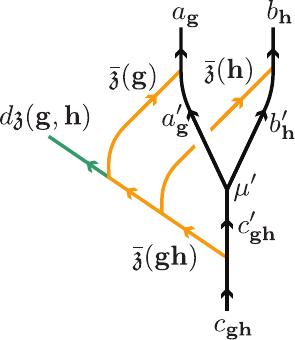}}}}}

\newcommand{\TorsorCompOne}{\mathord{\vcenter{\hbox{\includegraphics[scale=1]{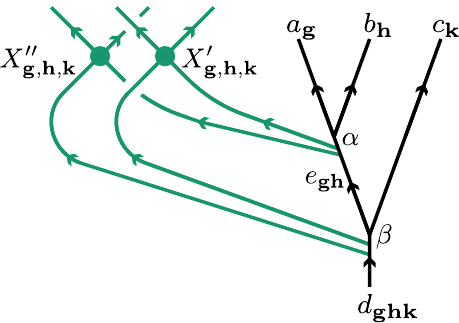}}}}}
\newcommand{\TorsorCompTwo}{\mathord{\vcenter{\hbox{\includegraphics[scale=1]{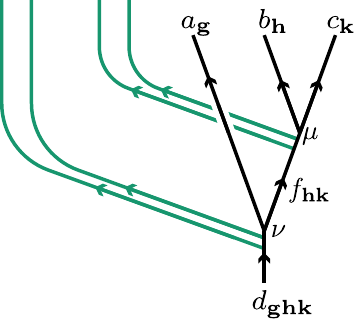}}}}}
\newcommand{\TorsorCompOneFused}{\mathord{\vcenter{\hbox{\includegraphics[scale=1]{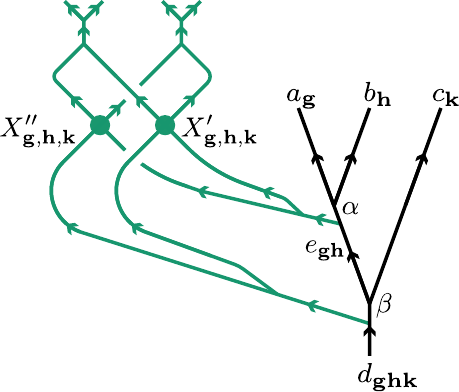}}}}}
\newcommand{\XFused}{\mathord{\vcenter{\hbox{\includegraphics[scale=1]{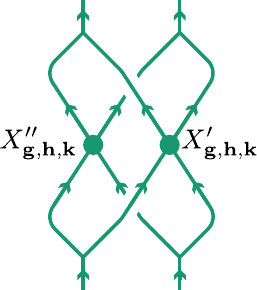}}}}}

\newcommand{\XFuseRHS}{\mathord{\vcenter{\hbox{\includegraphics[scale=1]{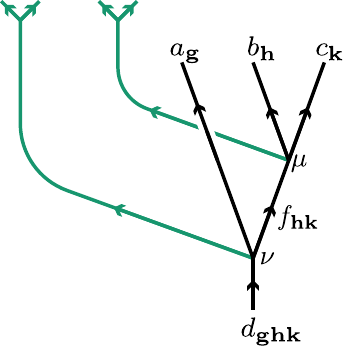}}}}}

\newcommand{\SymmAcGL}{\mathord{\vcenter{\hbox{\includegraphics[scale=1]{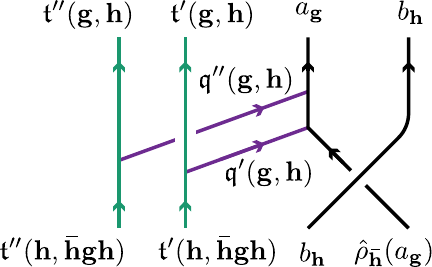}}}}}
\newcommand{\SymmAcGR}{\mathord{\vcenter{\hbox{\includegraphics[scale=1]{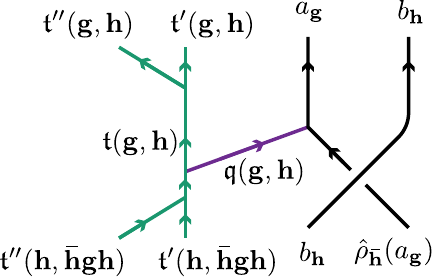}}}}}

\newcommand{\CommL}{\mathord{\vcenter{\hbox{\includegraphics[scale=1]{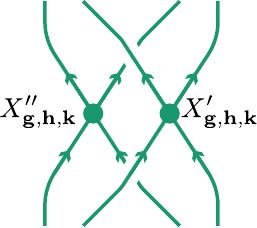}}}}}
\newcommand{\CommR}{\mathord{\vcenter{\hbox{\includegraphics[scale=1]{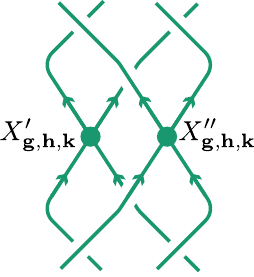}}}}}

\begin{document}

\title{
Torsorial actions on $G$-crossed braided tensor categories
}

\author{David Aasen}
\affiliation{Microsoft Station Q, Santa Barbara, California 93106-6105 USA}
\affiliation{Kavli Institute for Theoretical Physics, University of California, Santa Barbara, California 93106, USA}

\author{Parsa Bonderson}
\affiliation{Microsoft Station Q, Santa Barbara, California 93106-6105 USA}

\author{Christina Knapp}
\affiliation{Microsoft Station Q, Santa Barbara, California 93106-6105 USA}
\affiliation{Department of Physics and Institute for Quantum Information and Matter, California Institute of Technology, Pasadena, CA 91125, USA}
\affiliation{Walter Burke Institute for Theoretical Physics, California Institute of Technology, Pasadena, CA 91125, USA}

\date{\today}

\begin{abstract}
We develop a method for generating the complete set of basic data under the torsorial actions of $H^2_{[\rho]}(G,\mathcal{A})$ and $H^3(G,\text{U}(1))$ on a $G$-crossed braided tensor category $\mathcal{C}_G^\times$, where $\mathcal{A}$ is the set of invertible simple objects in the braided tensor category $\mathcal{C}$.
When $\mathcal{C}$ is a modular tensor category, the $H^2_{[\rho]}(G,\mathcal{A})$ and $H^3(G,\text{U}(1))$ torsorial action gives a complete generation of possible $G$-crossed extensions, and hence provides a classification.
This torsorial classification can be (partially) collapsed by relabeling equivalences that appear when computing the set of $G$-crossed braided extensions of $\mathcal{C}$.
The torsor method presented here reduces these redundancies by systematizing relabelings by $\mathcal{A}$-valued $1$-cochains.
We also use our methods to compute the composition rule of these torsor functors.
\end{abstract}

\maketitle
\tableofcontents

\section{Introduction}
\label{sec:intro}

$G$-crossed braided tensor categories generalize braided tensor categories to incorporate symmetry~\cite{turaev2000,Turaev2010,Etingof2010,Bark2019}.
In the presence of a symmetry action of the group $G$, a braided tensor category (BTC) $\bMTC$ is extended to a $G$-crossed braided tensor category $\bMTC_G^\times$, provided certain obstructions vanish.
When $\bMTC$ is a modular tensor category (MTC), $G$-crossed extensions are of particular relevance to the physics community for their characterization and classification of $(2+1)$D symmetry protected and symmetry enriched topological phases~\cite{Bark2019}.
The simple objects $a\in \bMTC$ correspond to the topologically distinct quasiparticles of the topological phase, while the simple objects $a_{\bf g} \in \bMTC_G^\times$ corresponding to nontrivial ${\bf g} \in G$ are associated with symmetry defects.

Compatibility of the topological data with the group action imposes a set of $G$-crossed consistency conditions, known as the pentagon and heptagon equations.
Classifying $G$-crossed BTCs amounts to solving these consistency conditions, which quickly becomes intractable for general $\bMTC$ and $G$.
Fortunately, the classification problem is substantially simplified through a connection to group cohomology~\cite{Etingof2010}.
For a given MTC and symmetry action $[\rho]$, there is a two-stage torsorial classification by $H^2_{[\rho]}(G,\mathcal{A})$ and $H^3(G,\text{U}(1))$, where the former takes coefficients in the set of Abelian simple objects in the MTC, which forms a group $\mathcal{A}$ under fusion.
While the $H^3(G,\text{U}(1))$ part of the torsorial action has been formulated concretely for the skeletonization of $G$-crossed BTCs in terms of gluing in SPT phases~\cite{Bark2019}, the analogous description of the torsorial action of $H^2_{[\rho]}(G,\mathcal{A})$ has been incomplete.
Moreover, the torsorial classification by cohomology groups sometimes overcounts distinct theories, as na\"ively different $G$-crossed BTCs can be identified through relabeling of the simple objects and gauge transformations.
This relabeling redundancy has previously been treated in an ad hoc manner for the $G$-crossed formalism.

In this work, we develop a method to describe the torsorial action of $Z^2_{[\rho]}(G,\mathcal{A})$ on the topological data of a $G$-crossed BTC $\bMTC_G^\times$.
In the special case where $\bMTC$ is a MTC, this torsor method generates the complete data of all distinct $G$-crossed extensions, given that of at least one such extension and the distinct $2$- and $3$-cocycles representing elements of $H^2_{[\rho]}(G,\mathcal{A})$ and $H^3(G,\text{U}(1))$.
The key insight is that torsoring in a $2$-cocycle $\mathfrak{t}\in Z^2_{[\rho]}(G, \mathcal{A})$ can modify the fusion ring of the theory.
By tracking the effect of this modification on fusion and braiding, we generate explicit expressions for the full topological data of the tosored theory.
Inserting the topological data into the $G$-crossed consistency conditions derives a relative obstruction $[\defectO_{r}]\in  H^4(G,\text{U}(1))$ that must vanish in order for the torsored theory to be fusion and $G$-crossed braiding consistent.
Additionally, we discuss the equivalences of na\"ively distinct $G$-crossed BTCs obtained by relabeling the simple objects of the theory.
By incorporating certain relabelings into the torsor method, we systematize their description in $G$-crossed classification.
Such relabeling equivalences can (partially) collapse the torsor structure of $G$-crossed BTCs, that is different elements of $H^2_{[\rho]}(G,\mathcal{A})$ and $H^3(G,\text{U}(1))$ do not necessarily correspond to distinct $G$-crossed BTCs when relabeling equivalences are included.
We also explicitly compute the compute the composition rule that allows the application of two torsor functors to be expressed as a single torsor functor, up to a gauge transformation.
Finally, we extend previous results~\cite{Bark2019} by using our torsor method to fully detail the topological data of all $G$-crossed extensions of MTCs upon which the symmetry acts trivially.

Our work builds on previous studies of the torsorial action of $Z^2_{[\rho]}(G,\mathcal{A})$.
The relative obstruction and torsorial action on the $F$-symbols were discussed in Refs.~\onlinecite{Etingof2010,Cui2016,Bark2019b}, while Ref.~\onlinecite{delaney2020} derived the full set of data for BTCs.
In contrast, the torsor method presented here generates the full set of topological data for $G$-crossed BTCs.

The remainder of this work is organized as follows.
In Section~\ref{sec:torsor_data}, we present our main result, the torsor method and the explicit basic data it generates.
In Section~\ref{sec:PW}, we analyze the consistency of this data and the corresponding obstructions.
In Section~\ref{sec:relabeling}, we consider equivalences of na\"ively different $G$-crossed BTCs and systematize certain relabeling equivalences by connecting them to the torsor method.
In Section~\ref{sec:Composition}, we give explicit expressions for the composition of two torsor functors.
Finally, in Section~\ref{sec:TrivSymmAc}, we provide an example use of the torsor method on theories with trivial symmetry action on the BTC.
Throughout this work, we employ the diagrammatic formalism and conventions for $G$-crossed BTCs used in Ref.~\onlinecite{Bark2019}, reviewed in Appendix~\ref{sec:background}.

\section{Basic data generated by the torsor method}
\label{sec:torsor_data}

The torsor method generates the topological data of a $G$-crossed BTC $\widehat{\bMTC}_G^\times$ given the topological data of another $G$-crossed extension $\bMTC_G^\times$ with the same symmetry action $[\rho]$, and the cochains relating the theories.
Abstractly, the torsor method provides a map between categories, described through a collection of functors
\begin{align}
\mathcal{F}_{\mathfrak{t},\Xdot}: \bMTC_G^\times \rightarrow \widehat{\bMTC}_G^\times,
\end{align}
parameterized by a representative 2-cochain $\mathfrak{t} \in Z^2_{[\rho]}(G,\mathcal{A})$ and a 3-cochain $\Xdot \in C^3(G,\text{U}(1))$.
When the trivial sector $\bMTC_{\bf 0}$ is an MTC, the torsor method translates the torsorial classification of $G$-crossed BTCs by $H^2_{[\rho]}(G,\mathcal{A})$ and $H^3(G,\text{U}(1))$ into explicit expressions of the torsorial action on the topological data.

We denote the quantities in $\widehat{\bMTC}_G^\times$ following the action of $\mathcal{F}_{\mathfrak{t},\Xdot}$ on ${\bMTC}_G^\times$ as $\widehat{Q} = \mathcal{F}_{\mathfrak{t},\Xdot}(Q)$, with the exception of the topological charge labels, which are unchanged by the torsor action.
Torsoring in the $2$-cocycle $\mathfrak{t}$ modifies the defect fusion rules as
\begin{align}
\label{torsorfusion}
a_{\bf g} \widehat{\otimes} b_{\bf h} \cong \mathfrak{t}({\bf g}, {\bf h}) \otimes (a_{\bf g} {\otimes} b_{\bf h}) = \bigoplus_{c_{\bf gh}} N_{a_{\bf g} b_{\bf h}}^{c_{\bf gh}} \mathfrak{t}({\bf g}, {\bf h}) \otimes c_{\bf gh} =  \bigoplus_{c_{\bf gh}} N_{a_{\bf g} b_{\bf h}}^{\bar{\mathfrak{t}}({\bf g}, {\bf h}) \otimes c_{\bf gh}} c_{\bf gh} .
\end{align}
In other words,
\begin{align}
\label{torsorfusioncoeff}
\widehat{N}_{a_{\bf g} b_{\bf h}}^{c_{\bf gh}} = N_{a_{\bf g} b_{\bf h}}^{\bar{\mathfrak{t}}({\bf g}, {\bf h}) \otimes c_{\bf gh}} .
\end{align}
Since $\mathfrak{t}({\bf g}, {\bf h})$ is always in $\mathcal{C}_{\bf 0}$, we will leave the subscript ${\bf 0}$ implicit; similarly, we will leave the subscript ${\bf 0}$ implicit when we write other $\mathcal{A}$-valued cochains.
It is worth observing that the new fusion rules are indeed associative when $\mathfrak{t} \in Z^{2}_{[\rho]}(G,\mathcal{A})$ (this condition is sufficient, but not necessary to achieve associativity).
The torsored fusion rules imply the fusion/splitting spaces of $\widehat{\bMTC}_{G}^\times$ are similarly modified.
In terms of the fusion spaces of the pre-torsored theory, we have
\begin{align}
\widehat{V}^{a_{\bf g} b_{\bf h}}_{c_{\bf gh}} \cong {V}^{\mathfrak{t}({\bf g},{\bf h}),a_{\bf g}, b_{\bf h}}_{c_{\bf gh}}
\cong {V}^{a_{\bf g} b_{\bf h}}_{\bar{\mathfrak{t}}({\bf g},{\bf h}) \otimes c_{\bf gh}}
.
\end{align}
Other choices for the isomorphisms between fusion spaces of $\widehat{\bMTC}_{G}^\times$ and ${\bMTC}_{G}^\times$, differing by $\mathfrak{t}({\bf g},{\bf h})$ would have worked equally well and can be related to the current choice by a braid; the choice here is most natural for the following way of defining the functor.

We define an explicit isomorphism between $\widehat{V}^{a_{\bf g} b_{\bf h}}_{c_{\bf gh}}$  and ${V}^{a_{\bf g} b_{\bf h}}_{c_{\bf gh}}$ diagrammatically by~\cite{Bark2019b} \begin{align}
\label{torsorfusionspace}
\FusionSpaceHat \cong \FusionSpaceR \equiv \FusionSpaceL
,
\end{align}
and introduce a shorthand (the diagram on the far right hand side) to simplify the diagrams.
The fusion diagrams on the right of the $\cong$ are defined according to the pre-torsored theory $\bMTC_G^\times$, and we indicate the ``phantom'' lines that disappear when passing to the left-hand side of the expression using green lines.
In the following, we will often leave the $\mathfrak{t}({\bf g},{\bf h})$ label implicit, since it can be inferred from the ${\bf g}$ and ${\bf h}$ labels emanating from the fusion vertex.
The dual vector spaces are found by reflecting the diagrams about the horizontal axis while preserving the orientation of the strands,
\begin{align}
\FusionSpaceHatdual \cong \FusionSpaceRDual \equiv \FusionSpaceLDual.
\end{align}

The $F$-moves of $\widehat{\bMTC}_G^\times$ are given by
\begin{align}
\FLeftHat
=\sum_{f_{\bf hk},\mu ,\nu} \left[ \widehat{F}^{a_{\bf g} b_{\bf h} c_{\bf k} }_{d_{\bf ghk}} \right]_{(e_{\bf gh},\alpha,\beta)(f_{\bf hk},\mu,\nu)}
\FRightHat
,
\end{align}
in which the $F$-symbols should be determined from the ${\bMTC}_G^\times$ data using the isomorphism of Eq.~(\ref{torsorfusionspace}).
The sum runs over all charges with the specified group label; the nontrivial elements are those obeying the new fusion rules.
However, in applying this isomorphism, the green $\mathfrak{t}$ lines require an additional operation for them to match correctly on both sides of the equation (in Ref.~\onlinecite{Bark2019b}, this was done by fusing them all to a single line).
For this, we write
\begin{align}
\label{eq:F-hat}
\!\!\!\!\FMoveCocyL\!\!\!\!\!\!
 =\sum_{f_{\bf hk},\mu ,\nu} \left[ \widehat{F}^{a_{\bf g} b_{\bf h} c_{\bf k} }_{d_{\bf ghk}} \right]_{(e_{\bf gh},\alpha,\beta)(f_{\bf hk},\mu,\nu)}
\FMoveCocyR
,
\end{align}
where the extra operator, which we call the ``cocycleator,'' is an isomorphism
\begin{align}
X_{{\bf g},{\bf h}, {\bf k}} : V^{\mathfrak{t}({\bf gh}, {\bf k}) , \mathfrak{t}({\bf g}, {\bf h})}_{\mathfrak{t}({\bf gh}, {\bf k}) \otimes \mathfrak{t}({\bf g}, {\bf h})} \rightarrow V^{\mathfrak{t}({\bf g}, {\bf hk}) , {}^{\bf g}\mathfrak{t}({\bf h}, {\bf k})}_{\mathfrak{t}({\bf g}, {\bf hk}) \otimes {}^{\bf g}\mathfrak{t}({\bf h}, {\bf k})}
\end{align}
The 2-cocycle condition ${}^{\bf g}\mathfrak{t}({\bf h}, {\bf k}) \bar{\mathfrak{t}}({\bf gh}, {\bf k}) \mathfrak{t}({\bf g}, {\bf hk}) \bar{\mathfrak{t}}({\bf g}, {\bf h}) = \I$  guarantees that such an isomorphism always exists.
Diagrammatically, we can define it as
\begin{equation}
\label{cocyleator}
\cocyleatorL\!\!\!\!X_{{\bf g},{\bf h},{\bf k}}~~ \equiv \Xdot ( {\bf g},{\bf h},{\bf k} ) \cocyleatorR ,
\end{equation}
with complex coefficients $\Xdot ( {\bf g},{\bf h},{\bf k} ) \in C^{3}(G,\text{U}(1))$ that will be constrained by the pentagon equation.
(We note that the group labels on $X$ determine all topological charge labels of the lines entering it, which can therefore be left implicit.)

We can explicitly write the $\widehat{F}$-symbols of $\widehat{\bMTC}_{G}^\times$ in terms of the basic data of $\bMTC_G^\times$ by evaluating the diagrams in Eq.~(\ref{eq:F-hat}), which yields
\begin{align}
\label{torFexplicit}
& \left[\widehat{F}^{a_{\bf g} b_{\bf h} c_{\bf k}} _{d_{\bf ghk}}\right]_{(e_{\bf gh},\alpha,\beta)(f_{\bf hk},\mu,\nu)} =
\sum_{\beta',\nu', \nu''}
\;
\left[{F}^{\mathfrak{t}({\bf g},{\bf h}) e'_{\bf gh} c_{\bf k}}_{d'_{\bf ghk}}\right]_{(e_{\bf gh},\beta), (d''_{\bf ghk},\beta')}
\left[{F}^{a_{\bf g} b_{\bf h} c_{\bf k}} _{d''_{\bf ghk}}\right]_{(e'_{\bf gh},\alpha, \beta'), (f'_{\bf hk},\mu ,\nu'')}
\nonumber \\
& \qquad
\times
\left( F^{\mathfrak{t}({\bf gh},{\bf k} ) \mathfrak{t}({\bf g},{\bf h}) d''_{\bf ghk}}_{d_{\bf ghk}} \right)^{-1}
F^{\mathfrak{t}({\bf g},{\bf hk}) {}^{\bf g}\mathfrak{t}({\bf h},{\bf k} ) d''_{\bf ghk}}_{d_{\bf ghk}}
\left[\left( F^{{}^{\bf g}\mathfrak{t}({\bf h},{\bf k}) a_{\bf g} f'_{\bf hk}}_{d'''_{\bf ghk}}  \right)^{-1}\right]_{(d''_{\bf ghk},\nu''), (a'_{\bf g},\nu')}
\nonumber \\
&
\qquad \qquad \times
\left( R^{{}^{\bf g} \mathfrak{t}({\bf h},{\bf k}) a_{\bf g}}_{ a'_{\bf g} }\right)^{-1}
\left[ F^{a_{\bf g} \mathfrak{t}({\bf h},{\bf k}) f'_{\bf hk}}_{d'''_{\bf ghk}} \right]_{(a'_{\bf g},\nu'), (f_{\bf hk},\nu)}
\Xdot({\bf g},{\bf h}, {\bf k})
\end{align}
where $e'_{\bf gh} = \bar{\mathfrak{t}}({\bf g},{\bf h})  \otimes e_{\bf gh}$, $f'_{\bf hk} = \bar{\mathfrak{t}}({\bf h},{\bf k}) \otimes f_{\bf hk}$, $a_{\bf g}' = {\mathfrak{t}}({\bf h},{\bf k}) \otimes a_{\bf g} = {}^{\bf g} {\mathfrak{t}}({\bf h},{\bf k}) \otimes a_{\bf g}$, $d'_{\bf ghk} = \bar{\mathfrak{t}}({\bf gh},{\bf k}) \otimes d_{\bf ghk}$, $d''_{\bf ghk}= \bar{\mathfrak{t}}({\bf g},{\bf h}) \otimes \bar{\mathfrak{t}}({\bf gh},{\bf k}) \otimes d_{\bf ghk} = {}^{\bf g}\bar{\mathfrak{t}}({\bf h},{\bf k}) \otimes \bar{\mathfrak{t}}({\bf g},{\bf hk}) \otimes d_{\bf ghk}$, and $d'''_{\bf ghk}= \bar{\mathfrak{t}}({\bf g},{\bf hk}) \otimes d_{\bf ghk}$. (Vertex labels of one-dimensional fusion spaces are left implicit.)
This is the same expression found in Ref.~\onlinecite{Bark2019b} for the $\widehat{F}$-symbols, with the exception of the $\Xdot({\bf g},{\bf h}, {\bf k})$ factor, which will generally be needed for the $\widehat{F}$-symbols to satisfy the pentagon equation (not just up to the obstruction).

We note that there are many possible ways of evaluating the diagrams that will generate different expressions for $\widehat{F}$-symbols in terms of the $\bMTC_G^\times$ data, by using different sequences of $F$-, $R$-, $U$-, and $\eta$-moves to relate the diagrams on the two sides of Eq.~(\ref{eq:F-hat}).
However, all such expressions will be equivalent as long as the $\bMTC_G^\times$ data satisfies the consistency conditions (pentagon and heptagons), which such a definition of $\widehat{F}$ well-defined.
Actually, this statement is slightly too strong, as the full set of consistency conditions are not required for $\widehat{F}$ to be well-defined; rather, the $\bMTC_G^\times$ data needs to satisfy all of the consistency conditions, except for the pentagon equation involving four nontrivial defects, i.e. Eq.~(\ref{eq:pentagoneqn}) for $a_{\bf g}$, $b_{\bf h}$, $c_{\bf k}$, $d_{\bf l}$ with ${\bf g}, {\bf h}, {\bf k}, {\bf l} \neq {\bf 0}$.
This is because Eq.~(\ref{eq:F-hat}) only involves at most three nontrivial defects, so the pentagon equations involving four nontrivial defects are never utilized to equate different diagrammatic sequences relating the two sides.

Since the braiding operators change the order of defects, the corresponding value of the green $\mathfrak{t}$ lines must transform accordingly.
Therefore, Eq.~(\ref{torsorfusionspace}) indicates braiding operators have basis elements corresponding to
\begin{align}
\RHatResolve
\end{align}

In light of this, we write the isomorphism of braiding operators as
\begin{align}
\label{eq:Rhat}
\widehat{R}^{a_{\bf g} b_{\bf h}} \cong \RHatOp
,
\end{align}
\begin{align}
\left(\widehat{R}^{a_{\bf g} b_{\bf h}}\right)^{-1} \cong \RHatOpInv
,
\end{align}
where
\begin{equation}\label{eq:q-def}
\mathfrak{q}({\bf g},{\bf h} ) = \bar{\mathfrak{t}}({\bf g},{\bf h} )\otimes \mathfrak{t}({\bf h},{\bf \bar{h}gh} )
\end{equation}
is the necessary charge transfer for the appropriate transformation of the green $\mathfrak{t}$ line through the braiding operation.
By consistency, this implies the new symmetry action on topological charges becomes
\begin{equation}
\label{eq:hat_rho}
{}^{\widehat{\bf k}} a_{\bf g} = \widehat{\rho}_{\bf k} (a_{\bf g}) =   {}^{\bf k}\bar{\mathfrak{q}}({\bf g},{\bf \bar{k}} ) \otimes {\rho}_{\bf k} (a_{\bf g})
.
\end{equation}
Notice that $\widehat{\rho}_{\bf k} (a_{\bf 0}) = {\rho}_{\bf k} (a_{\bf 0})$: as expected, the torsor method changes the fractionalization class while leaving the symmetry action unchanged.
While there is no freedom in defining $\mathfrak{q}$ for a given $\mathfrak{t}$, one can equivalently write the symmetry action as $\widehat{\rho}_{\bf k} (a_{\bf g}) =   \mathfrak{q}({\bf kg \bar{k}},{\bf k} ) \otimes {\rho}_{\bf k} (a_{\bf g})$ using the cocycle condition and trivial fusion with Abelian zero modes.
The torsored symmetry action $\widehat{\rho}$ satisfies the required conditions for a symmetry action in a $G$-crossed theory, namely:
\begin{enumerate}
\item $\hat{\rho}_{\bf k}(a_{\bf g} \widehat{\otimes} b_{\bf h}) \cong  \hat{\rho}_{\bf k}(a_{\bf g}) \widehat{\otimes}   \hat{\rho}_{\bf k}( b_{\bf h})$
\item $\hat{\rho}_{\bf h} ( \hat{\rho}_{\bf k}(a_{\bf g})) \cong \hat{\rho}_{\bf hk}(a_{\bf g})$
\item $a_{\bf g} \widehat{\otimes} b_{\bf h} \cong \hat{\rho}_{\bf g}(b_{\bf h}) \widehat{\otimes} a_{\bf g} \cong b_{\bf h} \widehat{\otimes} \hat{\rho}_{\bf \bar{h}}(a_{\bf g})$
\end{enumerate}

The $\widehat{R}$-moves are given by
\begin{align}
\RLeftHat =\sum_\nu \left[\widehat{R}^{a_{\bf g} b_{\bf h}}_{c_{\bf gh}} \right]_{\mu \nu} \RRightHat
.
\end{align}
Applying Eq.~(\ref{eq:Rhat}) to a vertex, we have the diagrammatic definition of the $\hat{R}$-symbols
\begin{align}
\label{eq:Rhat-symbol}
\RHatL \;\;=\sum_\nu \left[\widehat{R}^{a_{\bf g} b_{\bf h}}_{c_{\bf gh}} \right]_{\mu \nu}
\!\!\!\! \RHatR
.
\end{align}

We can explicitly write the $\widehat{R}$-symbols of $\widehat{\bMTC}_{G}^\times$ in terms of the basic data of $\bMTC_G^\times$ by evaluating the diagrams in Eq.~(\ref{eq:Rhat-symbol}), which yields
\begin{align}\label{eq:Rhatp}
\left[\widehat{R}^{a_{\bf g} b_{\bf h}}_{c_{\bf gh}} \right]_{\mu \nu}
=\sum_{\nu'}
\left[R^{a'_{\bf g} b_{\bf h}}_{c'_{\bf gh}} \right]_{\mu \nu'}
F^{\mathfrak{t}({\bf g },{\bf h})\mathfrak{q}({\bf g},{\bf h}) c'_{\bf gh}}_{c_{\bf gh}}
\left[ \left(F^{\mathfrak{q}({\bf g},{\bf h}) a'_{\bf g} b_{\bf h}}_{c''_{\bf gh}}\right)^{-1} \right]_{(c'_{\bf gh}, \nu')(a_{\bf g},\nu) }
,
\end{align}
where $a'_{\bf g} = \bar{\mathfrak{q}}({\bf g} , {\bf h}) \otimes a_{\bf g}$, $c'_{\bf gh} = \bar{\mathfrak{t}}({\bf h},{\bf \bar{h}gh}) \otimes c_{\bf gh}$, and $c''_{\bf gh} = \bar{\mathfrak{t}}({\bf g},{\bf h}) \otimes c_{\bf gh}$.

The $\widehat{U}$- and $\widehat{\eta}$-symbols can, in principle, be computed in terms of the $\bMTC_G^\times$ data by inserting the $\widehat{F}$- and $\widehat{R}$-symbols from Eqs.~(\ref{torFexplicit}) and (\ref{eq:Rhat-symbol}) into the heptagon equations.
However, doing so yields complicated expressions and an over-determined system of equations, where different topological charge and vertex labels in the heptagon equations may correspond to the same $\widehat{U}$- or $\widehat{\eta}$-symbols.
Directly verifying that all the different resulting expressions are equivalent is excruciating.
We instead use the following method that adroitly factors out operations acting on the phantom (green and purple) lines, thereby finding relatively simpler expressions for $\widehat{U}$- or $\widehat{\eta}$-symbols, and simplifying verification of the heptagon equations.

We define the $\widehat{U}$-symbols diagrammatically as
\begin{align}
\label{eq:Uhat-def}
\UHatRPrime \!\!\!\!\!=\sum_\nu
\left[ \widehat{U}_{\bf k} ( a_{\bf g},  b_{\bf h} ; c_{\bf gh}) \right]_{\mu \nu}\!\!\!\!\!\!\!\!\!\!\!\!\!\!\!\!\!\!
 \UHatLPrime
,
\end{align}
where the $P^{U}_{{\bf g},{\bf h},{\bf k}}$ operator is defined to be
\begin{align}
\label{eq:Psolution}
\PUdef = {U}_{\bf g}\left( {}^{\bf g}\mathfrak{t}({\bf h},{\bf k}) , {}^{\bf g}\mathfrak{q}({\bf h},{\bf k}) \right)  \PUdefRHS
.
\end{align}
This definition of $P^{U}_{{\bf g},{\bf h},{\bf k}}$ is justified in Sec.~\ref{sec:PW}, where we show that it is the necessary and sufficient for the heptagon equations involving $\widehat{U}$ to be satisfied.

We can explicitly write the $\widehat{U}$-symbols of $\widehat{\bMTC}_{G}^\times$ in terms of the basic data of $\bMTC_G^\times$ by evaluating the diagrams in Eq.~\eqref{eq:Uhat-def} and \eqref{eq:Psolution}, which yields
\begin{align}\label{eq:Uhat}
& \left[\widehat{U}_{\bf k}\left( a_{\bf g}, b_{\bf h}; c_{\bf gh} \right)\right]_{\mu \nu} =
\sum_{\nu',\nu'',\nu'''} \left[{U}_{\bf k}\left( a'_{\bf g}, b'_{\bf h}; c'''_{\bf gh} \right)\right]_{\mu \nu'}
{U}_{\bf k}\left( {}^{\bf k}\mathfrak{t}({\bf \bar{k}gk},{\bf \bar{k}hk}) , c'''_{\bf gh} ; c'_{\bf gh} \right)
\notag \\
& \qquad \times
{U}_{\bf g}\left( {}^{\bf g}\mathfrak{t}({\bf h},{\bf k}) , {}^{\bf g}\mathfrak{q}({\bf h},{\bf k}) \right)
\left( R^{{}^{\bf g} \mathfrak{q}({\bf h},{\bf k}) a_{\bf g}} \right)^{-1}
\left[ \left(F^{\mathfrak{q}({\bf g},{\bf k}) a'_{\bf g} b'_{\bf h}}_{c''''_{\bf gh}}\right)^{-1} \right]_{(c'''_{\bf gh}, \nu')(a_{\bf g},\nu'') }
\notag \\
& \qquad \times
\left[ \left(F^{{}^{\bf g}\mathfrak{q}({\bf h},{\bf k}) a_{\bf g} b'_{\bf h}}_{c''_{\bf gh}}\right)^{-1} \right]_{(c''''_{\bf gh}, \nu'')([\mathfrak{q}({\bf h},{\bf k})a_{\bf g}],\nu''') }
\left[ F^{a_{\bf g} \mathfrak{q}({\bf h},{\bf k}) b'_{\bf h}}_{c''_{\bf gh}} \right]_{([\mathfrak{q}({\bf h},{\bf k})a_{\bf g}],\nu''')(b_{\bf h}, \nu) }
\notag \\
& \qquad \times
\frac{F^{[\mathfrak{t}({\bf gk},{\bf \bar{k}hk})  \mathfrak{t}({\bf g},{\bf k})] \mathfrak{q}({\bf g},{\bf k})  c'''_{\bf gh}  }
F^{\mathfrak{t}({\bf g},{\bf hk}) {}^{\bf g}\mathfrak{t}({\bf h},{\bf k}) {}^{\bf g}\mathfrak{q}({\bf h},{\bf k})  }
F^{[{}^{\bf g}\mathfrak{t}({\bf h},{\bf k}) \mathfrak{t}({\bf g},{\bf hk}) ] {}^{\bf g}\mathfrak{q}({\bf h},{\bf k})   c''''_{\bf gh}  }
F^{\mathfrak{t}({\bf gh},{\bf k}) \mathfrak{t}({\bf g},{\bf h})  c''_{\bf gh}  }
}
{F^{\mathfrak{t}({\bf k},{\bf \bar{k}ghk}) {}^{\bf k}\mathfrak{t}({\bf \bar{k}gk},{\bf \bar{k}hk})  c'''_{\bf gh}  }
F^{\mathfrak{t}({\bf gk},{\bf \bar{k}hk}) \mathfrak{t}({\bf g},{\bf k})  \mathfrak{q}({\bf g},{\bf k})  }
F^{\mathfrak{t}({\bf gh},{\bf k}) \mathfrak{q}({\bf gh},{\bf k})  c'_{\bf gh}  }
}
\notag \\
& \qquad \times
\frac{\Xdot({\bf g},{\bf k}, {\bf \bar{k}hk}) }{\Xdot({\bf g},{\bf h}, {\bf k}) \Xdot({\bf k}, {\bf \bar{k}gk}, {\bf \bar{k}hk})}
,
\end{align}
where $a'_{\bf g} = \bar{\mathfrak{q}}({\bf g},{\bf k}) \otimes a_{\bf g}$, $b'_{\bf h} = \bar{\mathfrak{q}}({\bf h},{\bf k}) \otimes b_{\bf h}$, $c'_{\bf gh} = \bar{\mathfrak{q}}({\bf gh},{\bf k}) \otimes c_{\bf gh}$, $c''_{\bf gh} = \bar{\mathfrak{t}}({\bf g},{\bf h}) \otimes c_{\bf gh}$, $c'''_{\bf gh} = {}^{\bf k}\bar{\mathfrak{t}}({\bf \bar{k}gk},{\bf \bar{k}hk}) \otimes \bar{\mathfrak{q}}({\bf gh},{\bf k}) \otimes c_{\bf gh} = \bar{\mathfrak{t}}({\bf g},{\bf h}) \otimes \bar{\mathfrak{q}}({\bf g},{\bf k}) \otimes {}^{\bf g}\bar{\mathfrak{q}}({\bf h},{\bf k})\otimes c_{\bf gh}$, and $c''''_{\bf gh} = \bar{\mathfrak{t}}({\bf g},{\bf h}) \otimes {}^{\bf g}\bar{\mathfrak{q}}({\bf h},{\bf k})\otimes c_{\bf gh}$.
We leave labels uniquely determined by fusion implicit.

Similarly, we define the $\widehat{\eta}$-symbols diagrammatically as
\begin{align}\label{eq:etahat-def}
\etaHatRPrime =\hat{\eta}_{x_{\bf k}} ({\bf g},{\bf h}) \!\!\!\!\!\! \!\!\!\! \etaHatLPrime
,
\end{align}
where the $P^{U}_{{\bf g},{\bf h},{\bf k}}$ operator is defined to be
\begin{align}
\label{eq:Peta}
 \PEtaDef = U_{\bf g}({}^{\bf g}\mathfrak{t}({\bf \bar{g}kg},{\bf h}), {}^{\bf g}\mathfrak{q}({\bf \bar{g} k g},{\bf h}))^{-1}\PEtaDefR
,
\end{align}
and will show that this definition is necessary and sufficient for the heptagon equations involving $\widehat{\eta}$ to be satisfied.

We can explicitly write the $\widehat{\eta}$-symbols of $\widehat{\bMTC}_{G}^\times$ in terms of the basic data of $\bMTC_G^\times$ by evaluating the diagrams in Eq.~\eqref{eq:etahat-def} and \eqref{eq:Peta}, which yields
\begin{align}
\label{eq:etahat}
& \widehat{\eta}_{x_{\bf k}} \left( {\bf g}, {\bf h} \right) = \eta_{x'_{\bf k}} \left( {\bf g}, {\bf h} \right)
\frac{{U}_{\bf g}\left( {}^{\bf g}\mathfrak{q}({\bf \bar{g}kg},{\bf h}) , x'_{\bf k} \right)}{{U}_{\bf g}\left( {}^{\bf g}\mathfrak{t}({\bf \bar{g}kg},{\bf h}) , {}^{\bf g}\mathfrak{q}({\bf \bar{g}kg},{\bf h}) \right)}
R^{{}^{\bf k}\mathfrak{t}({\bf g},{\bf h}) x_{\bf k}}R^{ x'_{\bf k} \mathfrak{t}({\bf g},{\bf h})}
\notag \\
& \qquad \times
\frac{F^{\mathfrak{t}({\bf kg},{\bf h}) \mathfrak{t}({\bf k},{\bf g}) \mathfrak{q}({\bf k},{\bf g})  }
F^{\mathfrak{t}({\bf g},{\bf \bar{g}kgh}) {}^{\bf g}\mathfrak{t}({\bf \bar{g}kg},{\bf h}) {}^{\bf g}\mathfrak{q}({\bf \bar{g}kg},{\bf h})  }
F^{\mathfrak{t}({\bf gh},{\bf \bar{h}\bar{g}kgh}) \mathfrak{t}({\bf g},{\bf h}) x'_{\bf k} }
F^{\mathfrak{t}({\bf k},{\bf gh}) \mathfrak{q}({\bf k},{\bf gh})  x'''_{\bf k}  }
}{F^{\mathfrak{t}({\bf k},{\bf gh}) {}^{\bf k}\mathfrak{t}({\bf g},{\bf h}) x_{\bf k} }
F^{[\mathfrak{t}({\bf kg},{\bf h})  \mathfrak{t}({\bf k},{\bf g})] \mathfrak{q}({\bf k},{\bf g})  x''_{\bf k}  }
F^{[\mathfrak{t}({\bf g},{\bf \bar{g}kgh}) {}^{\bf g}\mathfrak{t}({\bf \bar{g}kg},{\bf h}) ] {}^{\bf g}\mathfrak{q}({\bf \bar{g}kg},{\bf h})   x'_{\bf k}  }
F^{\mathfrak{q}({\bf k},{\bf gh}) x'_{\bf k} \mathfrak{t}({\bf g},{\bf h}) }
}
\notag \\
& \qquad \times
\frac{\Xdot({\bf g},{\bf \bar{g}kg}, {\bf h}) }{\Xdot({\bf g},{\bf h}, {\bf \bar{h}\bar{g}kgh}) \Xdot({\bf k}, {\bf g}, {\bf h})}
,
\end{align}
where $x'_{\bf k} = \bar{\mathfrak{q}}({\bf k},{\bf gh}) \otimes x_{\bf k}$, $x''_{\bf k} = \bar{\mathfrak{q}}({\bf k},{\bf g}) \otimes x_{\bf k}$, and $x'''_{\bf k} = {\mathfrak{t}}({\bf g},{\bf h}) \otimes x'_{\bf k}$.

As was the case for the $\widehat{F}$-symbols, the definitions of the $\widehat{R}$-, $\widehat{U}$-, and $\widehat{\eta}$-symbols given above require that the $\bMTC_G^\times$ data satisfies the pentagon and heptagon consistency conditions involving no more than three nontrivial defects, in order to be well-defined (i.e. to not depend on a choice of evaluation path).
In particular, the pentagon equation involving four nontrivial defects is not required, as it is never utilized to relate different sequences of moves relating the two sides of Eqs.~(\ref{eq:Rhat-symbol}), (\ref{eq:Uhat-def}), and (\ref{eq:etahat-def}).

We note that the effect of the torsor action on $\bMTC$ should be trivial, except for changing the fractionalization class. In particular, we can verify that $\widehat{F}^{a_{\bf 0} b_{\bf 0} c_{\bf 0}}_{d_{\bf 0}} = {F}^{a_{\bf 0} b_{\bf 0} c_{\bf 0}}_{d_{\bf 0}}$, $\widehat{R}^{a_{\bf 0} b_{\bf 0}}_{ c_{\bf 0}}=R^{a_{\bf 0} b_{\bf 0}}_{ c_{\bf 0}}$, and $\widehat{U}_{\bf k}\left( a_{\bf 0}, b_{\bf 0}; c_{\bf 0} \right) = {U}_{\bf k}\left( a_{\bf 0}, b_{\bf 0}; c_{\bf 0} \right)$, while
\begin{align}
\label{eq:etahat_C0}
\widehat{\eta}_{x_{\bf 0}}({\bf g},{\bf h}) =  {\eta}_{x_{\bf 0}}({\bf g},{\bf h}) M_{x_{\bf 0} \mathfrak{t}({\bf g},{\bf h})}
,
\end{align}
where $M_{ab}= R^{ab}R^{ba}$ for Abelian topological charges $a$ and $b$, as expected.

\section{Consistency and relative obstruction}
\label{sec:PW}

The methods used in Section~\ref{sec:torsor_data} provide the data of $\widehat{\bMTC}_G^\times$ in terms of $\mathfrak{t} \in Z^2_{[\rho]}(G,\mathcal{A})$ and the data of $\bMTC_G^\times$, up to factors of the cocycleator $\Xdot$.
These factors were introduced to encapsulate the ability of the new data to satisfy the $G$-crossed consistency conditions.
Indeed, not all $\Xdot \in C^3(G,\text{U}(1))$ will lead to data that satisfy the full consistency conditions, and moreover it is not always possible to choose $\Xdot$ so that it does.
The obstruction to doing so is specifically manifested in this method as a failure of the $\widehat{F}$-symbols to satisfy the pentagon equation involving four nontrivial defects, while the heptagon equations are automatically satisfied, as we now show.

\begin{figure*}[t!]
\centering
\includegraphics[width=.99\textwidth]{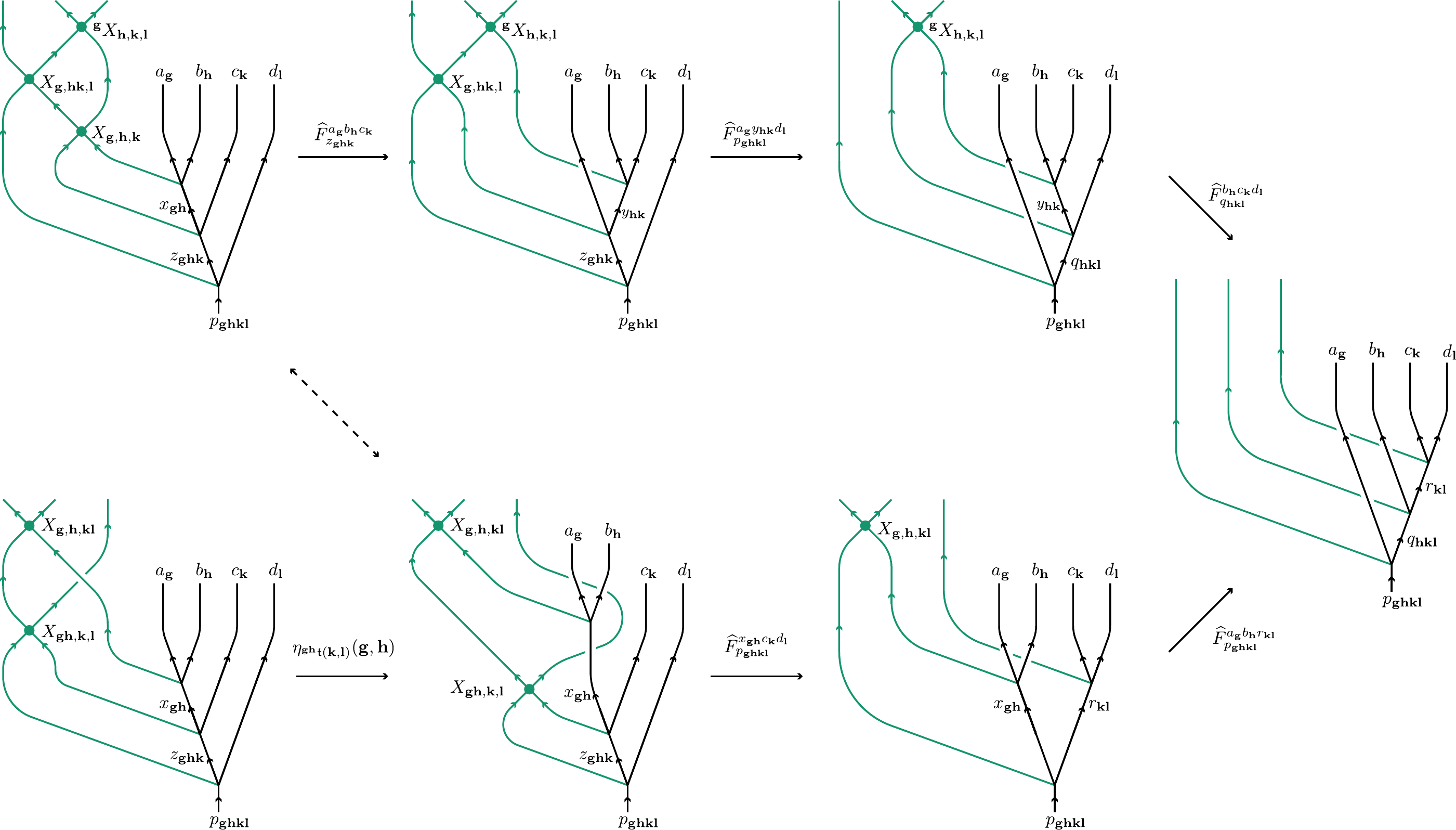}
\caption{
A sequence of diagrammatic moves in ${\bMTC}_G^\times$ involving the $\widehat{F}$-moves.
When the two diagrams connected by the dashed double-sided arrow are equal, the sequence becomes the pentagon equation for the $\widehat{F}$-moves in $\widehat{\bMTC}_G^\times$.
Equating these two diagrams is equivalent to requiring the condition that $\Xdot$ satisfies $\text{d} \Xdot  = \defectO_{r}^{-1}$.
}
\label{torsoredpentagon}
\end{figure*}

We first consider the sequence of diagrammatic moves shown in Figure~\ref{torsoredpentagon}.
By inspection of the figure, we see that if the ${\bMTC}_G^\times$ satisfies all the consistency conditions, then the $\widehat{F}$-symbols will satisfy the pentagon equation iff the cocyclator satisfies the condition
\begin{align}
\label{PWt}
\XConsisL \!\!\!\!\!=\frac{1}{\eta_{{}^{\bf gh}\mathfrak{t}({\bf k},{\bf l})}({\bf g},{\bf h})} \XConsisR
,
\end{align}
which equates the diagrams connected by the dashed double-sided arrow.
In the above, we write ${}^{\bf g}X_{{\bf h},{\bf k},{\bf l}}$ to mean that all strands emanating from the corresponding junction have been acted on by $\rho_{\bf g}$; the corresponding coefficient is ${}^{\bf g}\Xdot({\bf h},{\bf k},{\bf l})  \equiv  \frac{U_{\bf g}({}^{\bf g} \mathfrak{t}({\bf h},{\bf kl}),{}^{\bf gh}\mathfrak{t}({\bf k},{\bf l}))}{U_{\bf g}({}^{\bf g}\mathfrak{t}({\bf hk},{\bf l}),{}^{\bf g}\mathfrak{t}({\bf h},{\bf k}))}\Xdot({\bf h},{\bf k},{\bf l})$.

Rearranging Eq.~\eqref{PWt} so that all $\Xdot$-symbols appear on one side of the expression and all dependence on $\mathfrak{t}$ and the basic data of $\bMTC_G^\times$ appears on the other, we find
\begin{align}
\label{eq:X-pw}
\text{d} \Xdot =\defectO_{r}(\mathfrak{t})^{-1}
\end{align}
where $\text{d}$ is the usual coboundary operator, i.e., $\text{d}\Xdot({\bf g},{\bf h},{\bf k},{\bf l}) = \frac{\Xdot({\bf h},{\bf k},{\bf l})\Xdot({\bf g},{\bf hk},{\bf l})\Xdot({\bf g},{\bf h},{\bf k})} {\Xdot({\bf gh},{\bf k},{\bf l}) \Xdot({\bf g},{\bf h},{\bf kl})}$ and $\defectO_{r}(\mathfrak{t})$ is the relative obstruction.
A compact form of the relative obstruction first appeared in Ref.~\onlinecite{Etingof2010}, and was subsequently unpacked into an explicit formula in the mathematics~\cite{Cui2016} and physics~\cite{Bark2019b} literature, which we reproduce here by evaluating Eq.~(\ref{eq:X-pw}):
\begin{align}
\label{pwformula}
\defectO_{r}(\mathfrak{t})({\bf g},{\bf h},{\bf k},{\bf l})&= \eta_{{}^{\bf gh}\mathfrak{t}({\bf k},{\bf l})}({\bf g},{\bf h})
\frac{U_{\bf g}({}^{\bf g} \mathfrak{t}({\bf h},{\bf kl}),{}^{\bf gh}\mathfrak{t}({\bf k},{\bf l})) }
{U_{\bf g}({}^{\bf g}\mathfrak{t}({\bf hk},{\bf l}),{}^{\bf g}\mathfrak{t}({\bf h},{\bf k}))}
R^{{}^{\bf gh} \mathfrak{t}({\bf k},{\bf l}) \mathfrak{t}({\bf g},{\bf h})}
\nonumber \\
&\qquad \times
\frac{F^{\mathfrak{t}({\bf gh},{\bf kl}) \mathfrak{t}({\bf g},{\bf h}) {}^{\bf gh}\mathfrak{t}({\bf k},{\bf l})}}
{F^{\mathfrak{t}({\bf gh},{\bf kl})  {}^{\bf gh}\mathfrak{t}({\bf k},{\bf l})\mathfrak{t}({\bf g},{\bf h})}}
\frac{F^{\mathfrak{t}({\bf g},{\bf hkl}) {}^{\bf g}\mathfrak{t}({\bf hk},{\bf l} ) {}^{\bf g}\mathfrak{t}({\bf h},{\bf k})}}
{F^{\mathfrak{t}({\bf g},{\bf hkl}) {}^{\bf g} \mathfrak{t}({\bf h},{\bf kl}) {}^{\bf gh}\mathfrak{t}({\bf k},{\bf l}) }}
\frac{F^{\mathfrak{t}({\bf ghk},{\bf l}) \mathfrak{t}({\bf gh},{\bf k} ) \mathfrak{t}({\bf g},{\bf h})}}
{F^{\mathfrak{t}({\bf ghk},{\bf l}) \mathfrak{t}({\bf g},{\bf hk}) {}^{\bf g}\mathfrak{t}({\bf h},{\bf k} ) }}
.
\end{align}

There exist solutions for $\Xdot$ satisfying Eq.~\eqref{eq:X-pw} if and only if the relative obstruction is a coboundary $\defectO_{r}(\mathfrak{t}) \in B^4(G,\text{U}(1))$.
Clearly changing $\defectO_{r}(\mathfrak{t})$ by a coboundary does not change whether solutions exist.
Moreover, gauge transformations of the basic data generally change $\defectO_{r}(\mathfrak{t})$ by a coboundary.
In particular, vertex basis gauge transformations give
\begin{align}
\widetilde{\defectO}_{r}(\mathfrak{t}) &= \defectO_{r}(\mathfrak{t}) \text{d} \mu(\mathfrak{t}) , \\
\mu(\mathfrak{t}) ({\bf g},{\bf h},{\bf k}) & = \frac{\Gamma^{\mathfrak{t}({\bf g},{\bf hk}) \,^{\bf g}\mathfrak{t}({\bf h},{\bf k})} }{\Gamma^{\mathfrak{t}({\bf gh},{\bf k}) \mathfrak{t}({\bf g},{\bf h})  }}
,
\end{align}
while symmetry action gauge transformations give
\begin{align}
\check{\defectO}_{r}(\mathfrak{t}) &= \defectO_{r}(\mathfrak{t})  \text{d} \beta(\mathfrak{t}) , \\
\beta(\mathfrak{t}) ({\bf g},{\bf h},{\bf k}) & = \gamma_{\,^{\bf g}\mathfrak{t}({\bf h},{\bf k})}({\bf g})
,
\end{align}

Thus, it is natural to define the relative obstruction class as the equivalence class $[\defectO_{r}(\mathfrak{t})]$ related by coboundaries.
Moreover, it can be shown~\cite{Etingof2010} that the relative obstruction is a cocycle $\defectO_{r}(\mathfrak{t}) \in Z^4(G,\text{U}(1))$, so $[\defectO_{r}(\mathfrak{t})] \in H^4(G,\text{U}(1))$.

\begin{figure*}[t!]
\centering
\includegraphics[width=.99\textwidth]{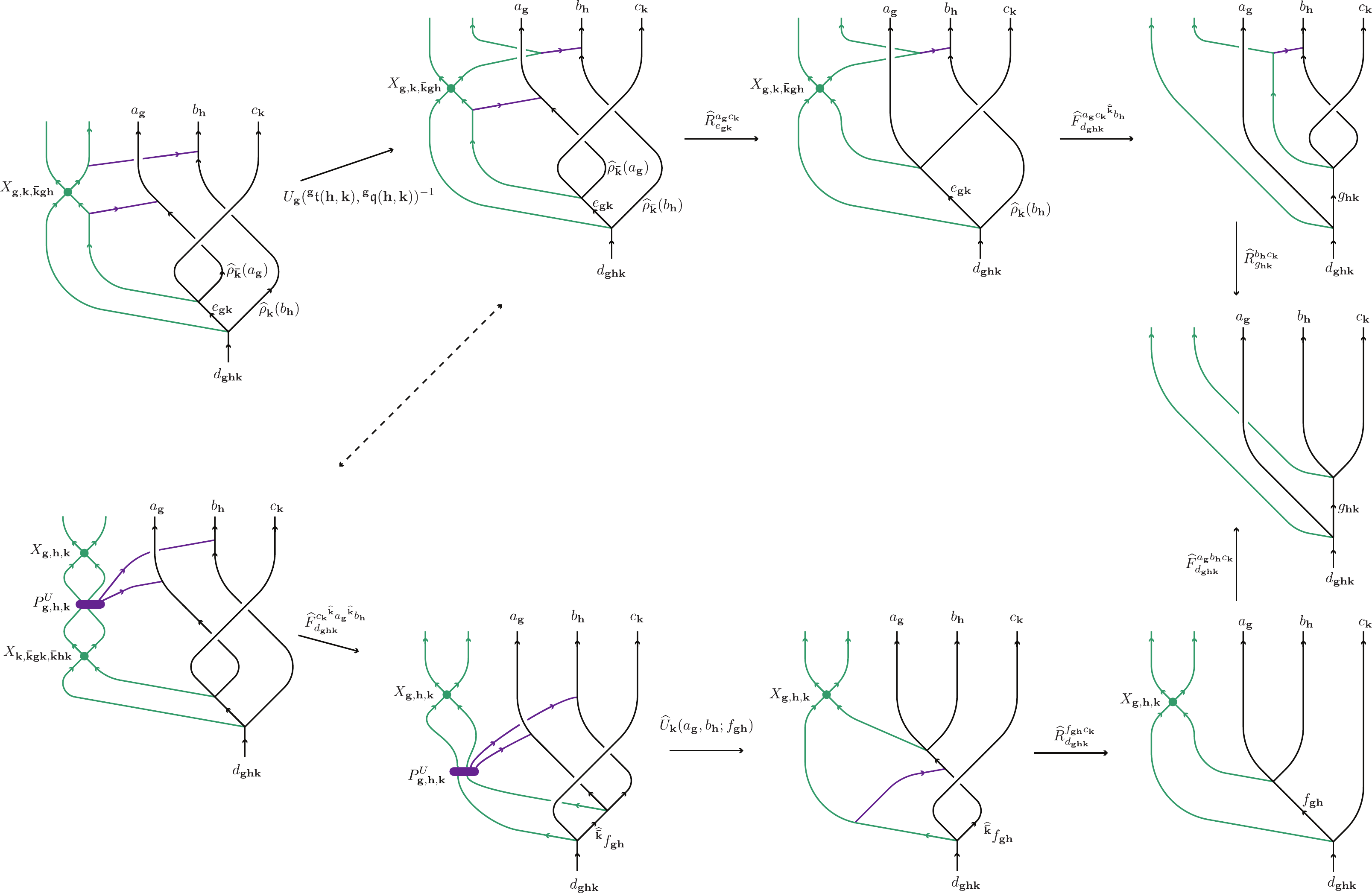}
\caption{A sequence of diagrammatic moves in ${\bMTC}_G^\times$ involving the $\widehat{F}$-, $\widehat{R}$-, and $\widehat{U}$-moves.
When the two diagrams connected by the dashed double-sided arrow are equal, the sequence becomes the $\widehat{U}$ heptagon equation for the $\widehat{\bMTC}_G^\times$ data.
Equating these two diagrams is equivalent to requiring the definition Eq.~(\ref{eq:Uhat-def}) for $P^{U}_{{\bf g},{\bf h},{\bf k}}$.
}
\label{fig:torsoredheptagonU}
\end{figure*}

\begin{figure*}[htp]
\centering
\includegraphics[width=.99\textwidth]{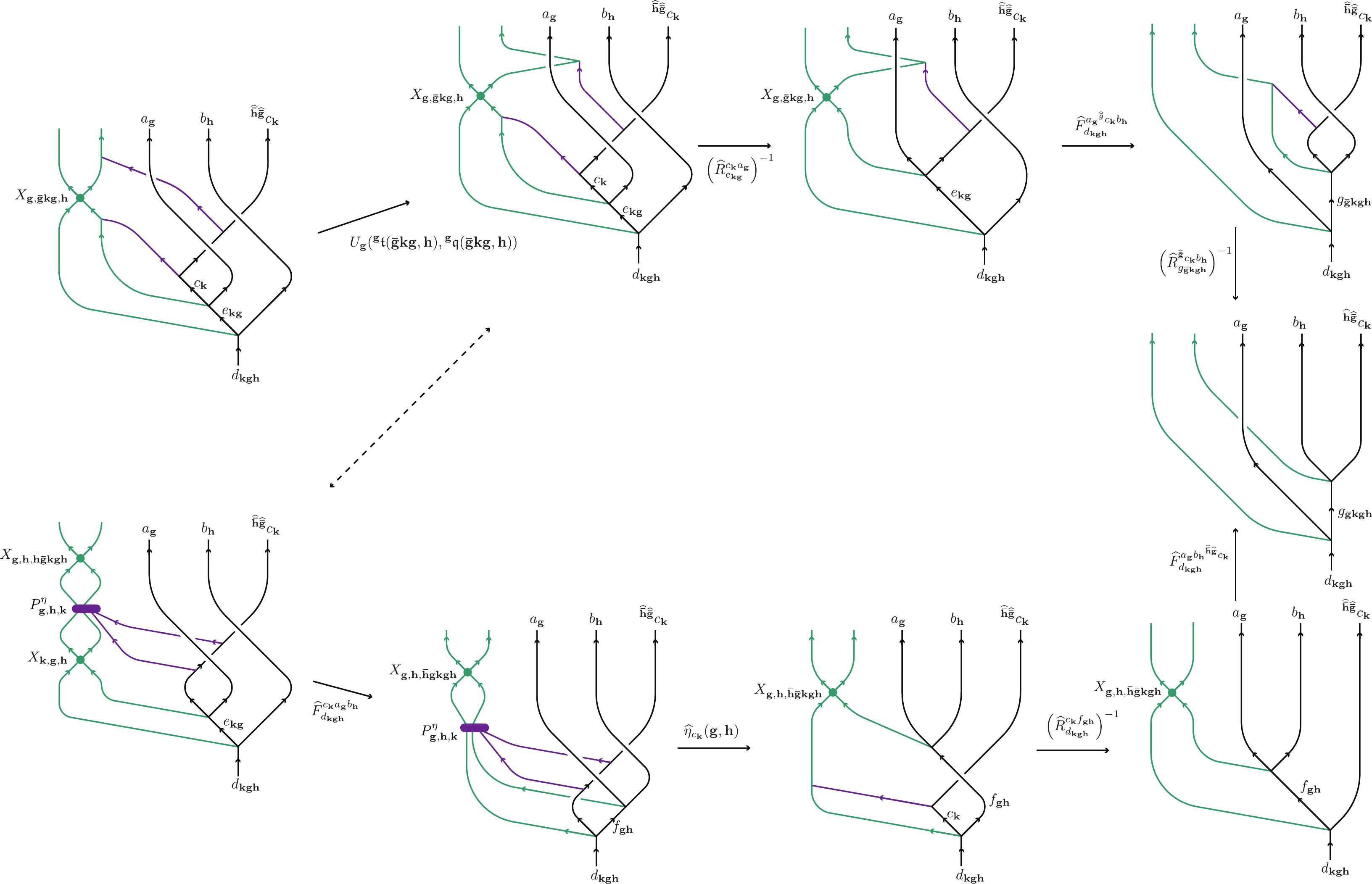}
\caption{A sequence of diagrammatic moves in ${\bMTC}_G^\times$ involving the $\widehat{F}$-, $\widehat{R}$-, and $\widehat{\eta}$-moves.
When the two diagrams connected by the dashed double-sided arrow are equal, the sequence becomes the $\widehat{\eta}$ heptagon equation for the $\widehat{\bMTC}_G^\times$ data.
Equating these two diagrams is equivalent to requiring the definition Eq.~(\ref{eq:etahat-def}) for $P^{\eta}_{{\bf g},{\bf h},{\bf k}}$.}
\label{fig:torsoredheptagoneta}
\end{figure*}

In order to show the data obtained for $\widehat{\bMTC}_G^\times$ satisfies the full consistency conditions, we must also verify the heptagon equations.
For the $U$-symbol heptagon equation, we consider the sequence of diagrammatic moves in Figure~\ref{fig:torsoredheptagonU}.
For the $\eta$-symbol heptagon equation, we consider the sequence of diagrammatic moves in Figure~\ref{fig:torsoredheptagoneta}.
By inspection of the figures, we see that if the ${\bMTC}_G^\times$ satisfies all the consistency conditions (for no more than three defects), then the data for $\widehat{\bMTC}_G^\times$ will satisfy the heptagon equations iff the $P^{U}_{{\bf g},{\bf h},{\bf k}}$ and $P^{\eta}_{{\bf g},{\bf h},{\bf k}}$ operators are defined as shown in Eqs.~(\ref{eq:Uhat-def}) and (\ref{eq:etahat-def}).
There is no obstruction that prevents these consistency conditions from being satisfied, assuming ${\bMTC}_G^\times$ satisfies all the consistency conditions for no more than three defects.

While the data cannot describe a consistent $(2+1)$D theory when the defectification obstruction class $[\defectO] \in H^4(G,\text{U}(1))$ is nontrivial, it is nonetheless meaningful to consider such data.
Physically, the data with nontrivial obstruction does not describe a strictly $(2+1)$D theory, but it can be used to describe an anomalous $(2+1)$D theory occurring at the boundary of a $(3+1)$D symmetry protected topological phase with bulk characterized by $[\defectO]$~\cite{Chen15}, the combination of which can be consistent due to the mechanism of anomaly inflow.
The obstruction class characterizes the anomaly of such theories.
Mathematically, we can see from the construction that the torsor action can be applied for $G$-crossed theories that have nontrivial defectification obstructions, as long as their topological data satisfies the heptagon equations and the pentagon equations for three or fewer nontrivial defects.
In this case, the torsor action also provides a transformation of the defectification obstruction of the theory.
Specifically, we can write this obstruction in terms of the $F$-symbols as
\begin{align}
\defectO({\bf g},{\bf h},{\bf k},{\bf l})&= F^{a_{\bf g} b_{\bf h} c_{\bf k} }_{z_{\bf ghk}} F^{a_{\bf g} y_{\bf hk} d_{\bf l} }_{p_{\bf ghkl}} F^{b_{\bf h} c_{\bf k} d_{\bf l} }_{q_{\bf ghkl}} \left(F^{a_{\bf g} b_{\bf h} r_{\bf kl} }_{p_{\bf ghkl}}\right)^{-1} \left(F^{x_{\bf gh} c_{\bf k} d_{\bf l} }_{p_{\bf ghkl}}\right)^{-1}
,
\end{align}
for the potentially obstructed/anomalous ${\bMTC}_G^\times$ theory.
Then the sequence of diagrammatic moves in Figure~\ref{torsoredpentagon} can be seen to lead to the defectification obstruction for $\widehat{\bMTC}_G^\times$ given by
\begin{align}
\widehat{\defectO}({\bf g},{\bf h},{\bf k},{\bf l}) &= \text{d}\Xdot({\bf g},{\bf h},{\bf k},{\bf l}) \defectO_{r}({\bf g},{\bf h},{\bf k},{\bf l}) \defectO({\bf g},{\bf h},{\bf k},{\bf l})
.
\end{align}
In terms of defectification obstruction classes (with equivalence under multiplication by coboundaries),
\begin{align}
[\widehat{\defectO}] &= [\defectO_{r} \defectO ]
,
\end{align}
which is why $\defectO_{r}$ is termed a relative obstruction.
When $\defectO_{r} \in B^4(G,\text{U}(1))$, $[\widehat{\defectO}] = [\defectO ]$, and it is possible to solve for $\Xdot$ such that the resulting data for $\widehat{\bMTC}_G^\times$ satisfies the full consistency conditions, including the pentagon equation with four nontrivial defects, if and only if the data for ${\bMTC}_G^\times$ does.

When an $\Xdot$ exists that yields an unobstructed $\widehat{\bMTC}_G^\times$, it is clear that $\Xdot' = \alpha \Xdot$ for any 3-cocycle $\alpha \in Z^{3}(G,\text{U}(1))$ will also yield an unobstructed theory, where $\alpha \Xdot$ denotes multiplication in $C^{3}(G,\text{U}(1))$, that is
\begin{align}\label{eq:X-alpha}
\alpha \Xdot ({\bf g},{\bf h},{\bf k}) = \alpha({\bf g},{\bf h},{\bf k}) \Xdot({\bf g},{\bf h},{\bf k}).
\end{align}
We recognize such shifts of $\Xdot$ as the action of the $Z^3(G,\text{U}(1))$ torsor corresponding to gluing in a bosonic SPT
\begin{align}
\widehat{\bMTC}_G^\times \mapsto  \spt_G^{[\alpha]} \underset{G}{\boxtimes}\widehat{\bMTC}_G^\times
,
\end{align}
where $\spt_G^{[\alpha]}$ are the $G$-crossed MTCs with trivial ${\bf 0}$-sector, as detailed in Appendix~\ref{sec:bSPT}.
In particular, the $3$-cocycle $\alpha$ enters the $\spt_G^{[\alpha]}$ data in Eqs.~(\ref{eq:F-alpha}-\ref{eq:eta-alpha}) in precisely the same way that $\Xdot$ enters the torsored $\widehat{\bMTC}_G^\times$ data in Section~\ref{sec:torsor_data}.
It follows that the torsor method presented here also incorporates the $Z^3(G,\text{U}(1))$ torsor action:
\begin{align}
\mathcal{F}_{\mathfrak{t}, \alpha \Xdot} = \spt_G^{[\alpha]} \underset{G}{\boxtimes}\mathcal{F}_{\mathfrak{t}, \Xdot}
.
\end{align}

\section{Equivalences of \texorpdfstring{$G$}{G}-crossed braided tensor categories}
\label{sec:relabeling}

Two $G$-crossed BTCs whose topological data are not identical are not necessarily distinct theories, as there are various notions of equivalences of theories. This includes equivalence under vertex basis gauge transformations and symmetry action gauge transformations~\cite{Bark2019}, see the review in Appendix~\ref{sec:background}.
We denote the theory obtained by the action of such gauge transformations as $\mathcal{G}_{\Gamma,\gamma}\left({\bMTC}_G^\times \right)$.

Another equivalence of categories is that of relabeling the simple objects; if one simply makes a different choice of the symbol used to denote an object, it is still the same theory.
While this statement is obvious, drawing the labels from the same set, so that relabeling is a permutation, and including gauge transformations can obfuscate the relation.
For a relabeling that permutes the simple objects of a category $\mathcal{C}$ as $a \mapsto \rel(a)$, the relabeled topological data is defined in the obvious way, that is
\begin{align}
\rel \left( Q(a, b,\ldots) \right) &= Q(\rel(a), \rel(b),\ldots)
,
\end{align}
for a quantity $Q$ in the category, and the resulting theory (complete set of basic data) is denoted $\rel\left(\mathcal{C}\right)$.
When $ \varphi = \mathcal{G} \circ \rel$ leaves the topological data invariant for some choice of gauge transformation $\mathcal{G}$, $\varphi$ is an auto-equivalence map, which is considered a topological symmetry.
When the relabeling does not map back to topological data that is gauge equivalent to the original data, it indicates an equivalence between na\"ively distinct theories $\mathcal{C}$ and $\rel\left(\mathcal{C}\right)$.
In the case of $G$-crossed BTCs, relabelings are required to leave the group labels fixed, that is $\rel(a_{\bf g}) = \rel(a)_{\bf g}$.
Physically, one can understand this as following from the extrinsic nature of the global symmetries of a system, i.e. they can be measured and fixed with respect to some external reference.
Mathematically, one can relax this condition to include relabelings that act as elements of $\text{Aut}(G)$ on the symmetry group labels~\cite{Cain2017}.

There is an important subset of the possible relabelings of a $G$-crossed BTC, which are the set of relabelings parameterized by normalized 1-cochains $\mathfrak{z} \in C^1(G,\mathcal{A})$, where $\mathcal{A}$ is the set of Abelian charges in $\mathcal{C}_{\bf 0}$.
Specifically, these relabelings are given by
\begin{align}
\rel_{\mathfrak{z}} \left( a_{\bf g} \right) &= \mathfrak{z}({\bf g}) \otimes a_{\bf g} = a'_{\bf g}
.
\end{align}
Such relabelings automatically preserve the $G$-grading and yield new data whose fusion rules differ from the original by the $2$-coboundary $\text{d}\mathfrak{z} \in B^2(G,\mathcal{A})$, that is
\begin{align}
\rel_{\mathfrak{z}}  \left( N_{a_{\bf g} b_{\bf h}}^{c_{\bf gh}} \right) &= N_{a'_{\bf g} b'_{\bf h}}^{ c'_{\bf gh}} = N_{ a_{\bf g}  b_{\bf h}}^{ \text{d} \mathfrak{z}({\bf g} ,{\bf h} ) \otimes c_{\bf gh}} .
\end{align}
As such, one might expect these relabelings to be equivalent (up to a gauge transformation) to the torsor functor $\mathcal{F}_{\mathfrak{t},\Xdot}$ for $\mathfrak{t} = \text{d} \mathfrak{z}$ and some choice of $\Xdot$.
We now show this is indeed always true.

We begin by identifying an isomorphism of the torsored fusion space to the relabeled fusion space through the diagrammatic relation
\begin{align}
\label{vertexbasisiso}
\relabelbosonicR  \equiv \sum_{\mu'} \left[ \Gamma^{a'_{\bf g} b'_{\bf h}}_{c'_{\bf gh}} \right]_{\mu \mu'} \relabelbosonicL
.
\end{align}
This equation serves as a diagrammatic definition of a vertex basis gauge transformation $\Gamma$, the detailed value of which can be expressed in terms of the topological data by evaluating the diagrams.
(This evaluation is straightforward, but the details are not enlightening, so we will not present them.)

Then we notice that we can write
\begin{align}
\label{F-rel}
\RelabelLeftF
= \sum_{f',\mu',\nu'} \frac{\Xdot({\bf g},{\bf h},{\bf k})}{\mathscr{Z}_{\mathfrak{z}}({\bf g},{\bf h},{\bf k})}
\left[ F^{a'_{\bf{g}} b'_{\bf{h}} c'_{\bf{k}}}_{d'_{\bf{ghk}}}\right]_{(e'_{\bf gh},\alpha',\beta')(f'_{\bf hk},\mu',\nu')}
\RelabelRightF
,
\end{align}
where $\mathscr{Z}_{\mathfrak{z}}({\bf g},{\bf h},{\bf k})$ can be determined by evaluating the diagrams.
Thus, if we set
\begin{align}
\Xdot({\bf g},{\bf h},{\bf k}) &= \mathscr{Z}_{\mathfrak{z}}({\bf g},{\bf h},{\bf k}) \nonumber \\
&= \eta_{{}^{\bf gh}\bar{\mathfrak{z}}({\bf k})}({\bf g},{\bf h})
U_{\bf g}({}^{\bf g}\text{d}\mathfrak{z}({\bf h},{\bf k}),{}^{\bf g}\bar{\mathfrak{z}}({\bf h}))
U_{\bf g}({}^{\bf g}[\mathfrak{z}({\bf k})\bar{\mathfrak{z}}({\bf hk})],{}^{\bf gh}\bar{\mathfrak{z}}({\bf k}))
R^{{}^{\bf g} \text{d}\mathfrak{z}({\bf h},{\bf k}) , \bar{\mathfrak{z}}({\bf g})}
\nonumber \\
&\times
\frac{F^{\text{d}\mathfrak{z}({\bf gh},{\bf k}), {}^{\bf g}[\mathfrak{z}({\bf h})\bar{\mathfrak{z}}({\bf gh})], {}^{\bf g}\bar{\mathfrak{z}}({\bf h})}
F^{\text{d}\mathfrak{z}({\bf gh},{\bf k}), \text{d}\mathfrak{z}({\bf g},{\bf h}), \bar{\mathfrak{z}}({\bf g})}
F^{\text{d}\mathfrak{z}({\bf g},{\bf hk}), \bar{\mathfrak{z}}({\bf g}), {}^{\bf g}\text{d}\mathfrak{z}({\bf h},{\bf k})}
}
{F^{\text{d}\mathfrak{z}({\bf g},{\bf hk}), {}^{\bf g}\text{d}\mathfrak{z}({\bf h},{\bf k}), \bar{\mathfrak{z}}({\bf g})}
F^{[{}^{\bf g}\mathfrak{z}({\bf hk})\bar{\mathfrak{z}}({\bf ghk})],{}^{\bf g}\text{d}\mathfrak{z}({\bf h},{\bf k}), {}^{\bf g}\bar{\mathfrak{z}}({\bf h})}
F^{[{}^{\bf g}\mathfrak{z}({\bf hk})\bar{\mathfrak{z}}({\bf ghk})],{}^{\bf g}[\mathfrak{z}({\bf k})\bar{\mathfrak{z}}({\bf hk})],{}^{\bf gh}\bar{\mathfrak{z}}({\bf k})}
}
,
\label{eq:boundaryX}
\end{align}
we find that
\begin{align}
& \left[\widehat{F}^{a_{\bf g} b_{\bf h} c_{\bf k}} _{d_{\bf ghk}}\right]_{(e_{\bf gh},\alpha,\beta)(f_{\bf hk},\mu,\nu)}
\nonumber \\
& \qquad =
\sum_{\alpha',\beta',\mu',\nu'}
\left[ \Gamma^{a' b'}_{e'} \right]_{\alpha \alpha'} \left[ \Gamma^{e' c'}_{d'} \right]_{\beta \beta'}
\left[ F^{a'_{\bf{g}} b'_{\bf{h}} c'_{\bf{k}}}_{d'_{\bf{ghk}}}\right]_{(e'_{\bf gh},\alpha',\beta')(f'_{\bf hk},\mu',\nu')}
\left[ \Gamma^{b' c'}_{f'} \right]^{-1}_{\mu' \mu} \left[ \Gamma^{a' f'}_{d'} \right]^{-1}_{\nu' \nu}
.
\end{align}
In other words, $\mathcal{F}_{\text{d}\mathfrak{z},\Xdot} (F) = \mathcal{G}_{\Gamma , \gamma} \circ \rel_{\mathfrak{z}} (F)$. (This holds for any $\gamma$, since symmetry action gauge transformations do not affect the $F$-symbols.)
Comparing with the torsor method, we see that $\text{d}\mathscr{Z}_{\mathfrak{z}} = \defectO_{r}^{-1}$.

Next, we notice that we can write
\begin{align}
\label{R-rel}
\RRelabelL= \sum_{\nu'} \mathcal{Y}_{\mathfrak{z}}({\bf g},{\bf h},a'_{\bf g}) \left[ R^{a'_{\bf{g}} b'_{\bf{h}}}_{c'_{\bf{gh}}}\right]_{\mu'\nu'}\RRelabelR
,
\end{align}
where $\mathcal{Y}_{\mathfrak{z}}({\bf g},{\bf h},a'_{\bf g})$ can be determined by evaluating the diagrams.
(This evaluation is also straightforward, but not enlightening.)
If we set
\begin{align}
\gamma_{a'_{\bf g}}({\bf h}) = \mathcal{Y}_{\mathfrak{z}}({\bf g},{\bf h},a'_{\bf g}),
\end{align}
we find that
\begin{align}
\left[\widehat{R}^{a_{\bf g} b_{\bf h}}_{c_{\bf gh}} \right]_{\mu \nu} = \sum_{\mu',\nu'} \gamma_{a'_{\bf g}}({\bf h}) \left[ \Gamma^{b' {}^{\bar{h}}a'}_{c'} \right]_{\mu \mu'}
\left[ R^{a'_{\bf{g}} b'_{\bf{h}}}_{c'_{\bf{gh}}}\right]_{\mu'\nu'}
\left[ \Gamma^{a' b'}_{c'} \right]^{-1}_{\nu' \nu}
\end{align}
 $\mathcal{F}_{\text{d}\mathfrak{z},\Xdot} (R) = \mathcal{G}_{\Gamma, \gamma} \circ \rel_{\mathfrak{z}} (R)$.

It follows that the corresponding $U$- and $\eta$-symbols of the torsored and relabeled theories are also related through this gauge transformation.
Thus, we have shown that
\begin{equation}
\label{eq:boundeqrel}
\mathcal{G}_{\Gamma, \gamma} \circ \rel_{\mathfrak{z}} ({\bMTC}_G^\times) = \mathcal{F}_{\text{d} \mathfrak{z},\mathscr{Z}_{\mathfrak{z}}} ({\bMTC}_G^\times)
\end{equation}
for $\mathscr{Z}_{\mathfrak{z}}$ given in Eq.~(\ref{eq:boundaryX}) and the gauge transformations defined above.
Since relabeling objects and gauge transformations do not change whether or not the data satisfies the consistency conditions, this implies that Eq.~(\ref{eq:X-pw}) is satisfied and $[\widehat{\defectO}] = [\defectO]$, i.e. the relative obstruction class $[\defectO_{r}(\text{d} \mathfrak{z})]$ is trivial when $\mathfrak{t}$ is a coboundary.

An implication of Eq.~(\ref{eq:boundeqrel}) is that two $G$-crossed BTCs related by torsoring by a 2-coboundary $\mathfrak{t} = \text{d}\mathfrak{z} \in B^2_{[\rho]}(G,\mathcal{A})$ are equivalent, up to a $Z^3(G,\text{U}(1))$ torsor.
As was observed in Ref.~\onlinecite{Bark2019}, two $G$-crossed BTCs related by torsoring by a 3-coboundary $\alpha = \text{d}\varepsilon \in B^{3}(G,\text{U}(1))$ (i.e. gluing an SPT) are equivalent under the gauge transformation $\mathcal{G}_{\Gamma, \gamma}$ with
\begin{align}
\left[ \Gamma^{a_{\bf g} b_{\bf h}}_{c_{\bf gh}} \right]_{\mu \mu'} &= \varepsilon ({\bf g},{\bf h}) \delta_{\mu \mu'}, \\
\gamma_{a_{\bf g}}({\bf h}) &= \frac{\varepsilon ({\bf g},{\bf h})}{\varepsilon ({\bf h},{\bf \bar{h}gh})}
.
\end{align}
Thus, for the purposes of generating distinct $G$-crossed BTCs, it is natural to factor out these coboundary actions from the full $Z^2_{[\rho]}(G,\mathcal{A})$ and $Z^3(G,\text{U}(1))$  torsor action, yielding a $H^2_{[\rho]}(G,\mathcal{A})$ and $H^3(G,\text{U}(1))$ torsorial generation of theories.

In the case of $G$-crossed MTCs, the $H^2_{[\rho]}(G,\mathcal{A})$ and $H^3(G,\text{U}(1))$ torsor action was shown to be complete~\cite{Etingof2010}, yielding a classification of $G$-crossed extensions of a given MTC.
However, factoring out coboundaries in this way does not account for all possible relabeling equivalences, and so generally can provide an overcounting in the $G$-crossed MTC classification.
It follows from completeness that any relabeling on a $G$-crossed MTC ${\bMTC}_G^\times$ that acts as an auto-equivalence on ${\bMTC}_{\bf 0}$ must be gauge equivalent to a torsor functor; that is, if $\rel_{\bf 0} \in \text{Aut}({\bMTC}_{\bf 0})$, then $\mathcal{G}_{\Gamma,\gamma} \circ \rel = \mathcal{F}_{\mathfrak{t},\Xdot}$ for some $\Gamma$, $\gamma$, $\mathfrak{t}$, and $\Xdot$.

When a relabeling acts completely trivially with respect to ${\bMTC}_{\bf 0}$ and fractionalization, i.e. $\mathfrak{t} = \I$, it can potentially still change the defectification class, giving an equivalence of na\"ively distinct defectification classes for the same fractionalization class.
Indeed, the above analysis shows that $\rel_{\mathfrak{z}}$ when $\mathfrak{z} \in Z^1_{[\rho]}(G,\mathcal{A})$ (i.e. $\text{d}\mathfrak{z} = \I$) is gauge equivalent to $\mathcal{F}_{\I,\Xdot}$ where
\begin{align}
\label{eq:sptcollapse}
\Xdot({\bf g},{\bf h},{\bf k}) &= \frac{ \eta_{{}^{\bf gh}\bar{\mathfrak{z}}({\bf k})}({\bf g},{\bf h})
U_{\bf g}({}^{\bf g}[\mathfrak{z}({\bf k})\bar{\mathfrak{z}}({\bf hk})],{}^{\bf gh}\bar{\mathfrak{z}}({\bf k}))
}
{F^{[{}^{\bf g}\mathfrak{z}({\bf hk})\bar{\mathfrak{z}}({\bf ghk})],{}^{\bf g}[\mathfrak{z}({\bf k})\bar{\mathfrak{z}}({\bf hk})],{}^{\bf gh}\bar{\mathfrak{z}}({\bf k})}
}
,
\end{align}
is a 3-cocycle that can potentially correspond to a nontrivial $H^3(G,\text{U}(1))$ class.
In other words, there may exist 1-cocycle relabelings that equate na\"ively distinct defectification classes of the same fractionalization class.

In order for a relabeling to equate na\"ively distinct fractionalization classes, it must act as a nontrivial autoequivalence on ${\bMTC}_{\bf 0}$.
In particular, it must permute objects in the set of Abelian quasiparticles $\mathcal{A}$.
This permutation on $\mathcal{A}$ by the relabeling correspondingly equates fractionalization classes that are related by such relabeling.
For example, when the symmetry action is trivial ($\rho_{\bf g} = \openone$ for all ${\bf g}$), in which case the fractionalization classes can be identified with the cohomology classes $[\mathfrak{w}] \in H^2(G,\mathcal{A})$ (not just torsorially), the relabeling equates fractionalization classes $[\mathfrak{w}]$ with $[\rel(\mathfrak{w})]$.
Clearly, such relabelings cannot be parameterized by $G$ 1-cochains, which act trivially on ${\bMTC}_{\bf 0}$.
However, if the theory has an additional $K$ grading on the ${\bMTC}_{\bf 0}$-sector, then a similar construction can be applied with relabelings that are parameterized by $K$ 1-cochains (or perhaps $H$ 1-cochains, where $H$ is a $K$ extension of $G$).

\section{Composition of torsor functors}
\label{sec:Composition}

Two torsor functors can be composed to realize another torsor functor, up to a gauge transformation.
In particular, we will show that
\begin{align}
\label{eq:torsor_composition}
\mathcal{F}_{\mathfrak{t}'', \Xdot''} \circ \mathcal{F}_{\mathfrak{t}', \Xdot'} = \mathcal{G}_{\Gamma, \gamma} \circ \mathcal{F}_{\mathfrak{t}, \Xdot}
,
\end{align}
where the resulting torsor functor is specified by
\begin{align}
\label{eq:tcomp}
\mathfrak{t}({\bf g,h}) &= \mathfrak{t}'({\bf g,h}) \otimes \mathfrak{t}''({\bf g,h}),
\\
\label{eq:Xdotcomp}
\Xdot({\bf g,h,k}) &= \Xdot'({\bf g,h,k}) \Xdot''({\bf g,h,k}) \mathscr{Y}_{\mathfrak{t}', \mathfrak{t}''}({\bf g,h,k})
,
\end{align}
with
\begin{align}
& \mathscr{Y}_{\mathfrak{t}', \mathfrak{t}''}({\bf g,h,k}) =
\frac{1}{U_{\bf g}({}^{\bf g}\mathfrak{t}''({\bf h,k}), {}^{\bf g}\mathfrak{t}'({\bf h,k} ))}
\frac{R^{\mathfrak{t}''({\bf g,h})\mathfrak{t}'({\bf gh,k}) } }
{ R^{{}^{\bf g}\mathfrak{t}''({\bf h,k} ) \mathfrak{t}'({\bf g,hk}) } }
\frac{F^{\mathfrak{t}''({\bf gh,k})\mathfrak{t}'({\bf gh,k}) \mathfrak{t}''({\bf g,h}) } }
{F^{[\mathfrak{t}''({\bf gh,k})\mathfrak{t}'({\bf gh,k}) ] \mathfrak{t}''({\bf g,h}) \mathfrak{t}'({\bf g,h})}  }
\notag \\
& \qquad \times
\frac{
F^{[\mathfrak{t}''({\bf g,h}) \mathfrak{t}''({\bf gh,k})] \mathfrak{t}'({\bf gh,k}) \mathfrak{t}'({\bf g,h}) }
F^{\mathfrak{t}''({\bf g,hk}) {}^{\bf g} \mathfrak{t}''({\bf h,k}) \mathfrak{t}'({\bf g,hk} )}
F^{[\mathfrak{t}''({\bf g,hk}) \mathfrak{t}'({\bf g,hk}) ] {}^{\bf g}\mathfrak{t}''({\bf h,k}) {}^{\bf g}\mathfrak{t}'({\bf h,k}) }
}
{
F^{\mathfrak{t}''({\bf gh,k}) \mathfrak{t}''({\bf g,h})\mathfrak{t}'({\bf gh,k})}
F^{[\mathfrak{t}''({\bf g,h} ) \mathfrak{t}''({\bf gh,k} ) ]\mathfrak{t}'({\bf g,hk}) {}^{\bf g}\mathfrak{t}'({\bf h,k}) }
F^{\mathfrak{t}''({\bf g,hk}) \mathfrak{t}'({\bf g,hk}) {}^{\bf g} \mathfrak{t}''({\bf h,k})}
}
.
\label{eq:Ydot}
\end{align}
The gauge transformation $\mathcal{G}_{\Gamma, \gamma}$ corresponds to the vertex basis gauge transformation specified by
\begin{align}
\label{eq:Gamma_comp}
\left[\Gamma^{a_{\bf g} b_{\bf h}}_{e_{\bf gh}} \right]_{\mu \mu'} = \left(F^{\mathfrak{t}''({\bf g,h})\mathfrak{t}'({\bf g,h}) [\overline{\mathfrak{t}}({\bf g,h}) e_{\bf gh}] }_{e_{\bf gh}} \right)^{-1} \delta_{\mu \mu'}
,
\end{align}
together with the symmetry action gauge transformation specified by
\begin{align}
\label{eq:gamma_comp}
\gamma_{a_{\bf g}}({\bf h}) &=\frac{1}{R^{\mathfrak{q}''({\bf g,h}) \mathfrak{t}'({\bf g,h})}}
\frac{
F^{\mathfrak{t}''({\bf g,h}) \mathfrak{q}''({\bf g,h}) \mathfrak{t}'({\bf g,h}) }
F^{\mathfrak{t}({\bf g,h}) \mathfrak{q}''({\bf g,h}) \mathfrak{q}'({\bf g,h}) }
}
{
F^{\mathfrak{t}''({\bf h, \bar{h}gh}) \mathfrak{t}'({\bf g,h})\mathfrak{q}'({\bf g,h}) }
F^{\mathfrak{t}''({\bf g,h}) \mathfrak{t}'({\bf g,h}) \mathfrak{q}''({\bf g,h})}
}
F^{\mathfrak{q}''({\bf g,h}) \mathfrak{q}'({\bf g,h}) [\overline{\mathfrak{q}}({\bf g,h}) a_{\bf g}] }_{a_{\bf g}}
,
\end{align}
where $ \mathfrak{q}({\bf g,h})= \mathfrak{q}'({\bf g,h}) \otimes \mathfrak{q}''({\bf g,h})$.

We now derive Eqs.~\eqref{eq:torsor_composition}-\eqref{eq:Gamma_comp}.
Eq.~\eqref{eq:tcomp} clearly follows from Eqs.\eqref{torsorfusion}-\eqref{torsorfusionspace}.
To obtain  Eqs.~\eqref{eq:Xdotcomp}-\eqref{eq:Gamma_comp}, we use Eq.~\eqref{eq:F-hat} to write the diagrammatic description of the composition of torsor functors acting on the $F$-symbol as
\begin{align}
\label{eq:TorsorCompEqn}
\TorsorCompOne\!\!\!\!
=\!\!\sum_{f_{\bf hk}, \mu, \nu} \left[ \mathcal{F}_{\mathfrak{t}'', \Xdot''} \circ \mathcal{F}_{\mathfrak{t}', \Xdot'}(F)^{a_{\bf g} b_{\bf h}c_{\bf k} }_{d_{\bf ghk}}\right]_{(e_{\bf gh},\alpha, \beta)(f_{\bf hk}, \mu, \nu)}\!\!\!\! \TorsorCompTwo
.
\end{align}
Here, we explicitly draw the phantom lines attaching to the fusion channel lines below the vertices for clarity [recall Eq.~\eqref{torsorfusionspace}].
Next, we use diagrammatic relations to rewrite both sides of Eq.~\eqref{eq:TorsorCompEqn} to relate them to our standardized form for a torsor functor, as in Eq.~\eqref{eq:F-hat}.
We first fuse the phantom lines together to obtain
\begin{align}
\label{eq:TorsorCompGaugeTrans}
\TorsorCompOne =
\Gamma^{a_{\bf g} b_{\bf h}}_{e_{\bf gh} } \Gamma^{e_{\bf gh} c_{\bf k}}_{d_{\bf ghk} } \TorsorCompOneFused,
\end{align}
where $\Gamma$ are the inverse $F$-symbols given by Eq.~\eqref{eq:Gamma_comp}.
Then we rewrite the diagram on the right hand side of Eq.~\eqref{eq:TorsorCompGaugeTrans}, so that it has the form of the left hand side of Eq.~\eqref{eq:F-hat}, with a single cocycleator.
For this, we use
\begin{align}
\label{Xfused}
\XFused = \mathscr{W}({\bf g,h,k} ) \cocyleatorR = \frac{\mathscr{W}({\bf g,h,k} )}{\Xdot ( {\bf g},{\bf h},{\bf k} ) } \cocyleatorL\!\!\!\!X_{{\bf g},{\bf h},{\bf k}}
,
\end{align}
where the coefficient $\mathscr{W}$ is obtained by evaluating the diagram on the left hand side, and the (as of yet undetermined) coefficient $\Xdot ( {\bf g},{\bf h},{\bf k} )$ in the numerator on the right hand side simply cancels with the coefficient of the cocycleator $X_{{\bf g},{\bf h},{\bf k}}$.
Performing the diagrammatic computation of $\mathscr{W}$, we find
\begin{align}
& \mathscr{W}({\bf g,h,k}) =
\frac{R^{\mathfrak{t}''({\bf g,h})\mathfrak{t}'({\bf gh,k}) } }
{ R^{{}^{\bf g}\mathfrak{t}''({\bf h,k} ) \mathfrak{t}'({\bf g,hk}) } }
\frac{
F^{\mathfrak{t}''({\bf gh,k})\mathfrak{t}'({\bf gh,k}) \mathfrak{t}''({\bf g,h}) }
F^{[\mathfrak{t}''({\bf g,h}) \mathfrak{t}''({\bf gh,k})] \mathfrak{t}'({\bf gh,k}) \mathfrak{t}'({\bf g,h}) }
}
{F^{[\mathfrak{t}''({\bf gh,k})\mathfrak{t}'({\bf gh,k}) ] \mathfrak{t}''({\bf g,h}) \mathfrak{t}'({\bf g,h})}
F^{\mathfrak{t}''({\bf gh,k}) \mathfrak{t}''({\bf g,h})\mathfrak{t}'({\bf gh,k})}
}
\notag \\
& \qquad \times
\frac{
F^{\mathfrak{t}''({\bf g,hk}) {}^{\bf g} \mathfrak{t}''({\bf h,k}) \mathfrak{t}'({\bf g,hk} )}
F^{[\mathfrak{t}''({\bf g,hk}) \mathfrak{t}'({\bf g,hk}) ] {}^{\bf g}\mathfrak{t}''({\bf h,k}) {}^{\bf g}\mathfrak{t}'({\bf h,k}) }
}
{
F^{[\mathfrak{t}''({\bf g,h} ) \mathfrak{t}''({\bf gh,k} ) ]\mathfrak{t}'({\bf g,hk}) {}^{\bf g}\mathfrak{t}'({\bf h,k}) }
F^{\mathfrak{t}''({\bf g,hk}) \mathfrak{t}'({\bf g,hk}) {}^{\bf g} \mathfrak{t}''({\bf h,k})}
}
\Xdot'({\bf g,h,k}) \Xdot''({\bf g,h,k})
,
\label{eq:W}
\end{align}

Together, this yields
\begin{align}
&\TorsorCompOne
=  \Gamma^{a_{\bf g} b_{\bf h}}_{e_{\bf gh} } \Gamma^{e_{\bf gh} c_{\bf k}}_{d_{\bf ghk} }
\frac{\mathscr{W}({\bf g,h,k} )}{\Xdot ( {\bf g},{\bf h},{\bf k} ) }\FMoveCocyleLBranched
\notag \\
&\quad = \Gamma^{a_{\bf g} b_{\bf h}}_{e_{\bf gh} } \Gamma^{e_{\bf gh} c_{\bf k}}_{d_{\bf ghk} }
\frac{\mathscr{W}({\bf g,h,k} )}{\Xdot ( {\bf g},{\bf h},{\bf k} ) } \sum_{f_{\bf hk}, \mu, \nu} \left[  \mathcal{F}_{\mathfrak{t}, \Xdot}(F)^{a_{\bf g} b_{\bf h}c_{\bf k} }_{d_{\bf ghk}}\right]_{(e_{\bf gh},\alpha, \beta)(f_{\bf hk}, \mu, \nu)}
\XFuseRHS
\notag\\
&\quad =
\frac{\mathscr{W}({\bf g,h,k} )}{\Xdot ( {\bf g},{\bf h},{\bf k} ) } \frac{1}{U_{\bf g}({}^{\bf g}\mathfrak{t}''({\bf h,k}){}^{\bf g}\mathfrak{t}'({\bf h,k}) )}
\notag \\
&\quad \qquad \times
\sum_{f_{\bf hk}, \mu, \nu}
\frac{\Gamma^{a_{\bf g} b_{\bf h}}_{e_{\bf gh} } \Gamma^{e_{\bf gh} c_{\bf k}}_{d_{\bf ghk} }}{\Gamma^{b_{\bf h} c_{\bf k}}_{f_{\bf hk} } \Gamma^{a_{\bf g} f_{\bf hk} }_{d_{\bf ghk} } }
 \left[  \mathcal{F}_{\mathfrak{t}, \Xdot}(F)^{a_{\bf g} b_{\bf h}c_{\bf k} }_{d_{\bf ghk}}\right]_{(e_{\bf gh},\alpha, \beta)(f_{\bf hk}, \mu, \nu)}
\TorsorCompTwo
.
\label{eq:Fcomp}
\end{align}
Equating this with the right hand side of Eq.~\eqref{eq:TorsorCompEqn}, we obtain Eq.~\eqref{eq:torsor_composition} by setting
\begin{align}
\Xdot ( {\bf g},{\bf h},{\bf k} ) = \frac{\mathscr{W}({\bf g,h,k} )}{U_{\bf g}({}^{\bf g}\mathfrak{t}''({\bf h,k}){}^{\bf g}\mathfrak{t}'({\bf h,k}) )}
,
\end{align}
which yields Eqs.~\eqref{eq:Xdotcomp} and \eqref{eq:Ydot}.
The vertex basis part of the gauge transformation can be read off of Eq.~\eqref{eq:Fcomp}, which gives Eq.~\eqref{eq:Gamma_comp}.

The symmetry action part of the gauge transformation can be obtained through a similar analysis of the composition of torsor functors acting on the $R$-symbols.
While the $R$-symbols change under both vertex basis and symmetry action gauge transformations, the braiding operators do not change under vertex basis gauge transformations.
As such, we can use the braiding operators to more directly isolate the symmetry action part of $\mathcal{G}_{\Gamma, \gamma}$, which acts on the braiding operators as
\begin{align}
\mathcal{G}_{\Gamma, \gamma} \left( R^{a_{\bf g} b_{\bf h}} \right) = \gamma_{a_{\bf g}}({\bf h} ) R^{a_{\bf g} b_{\bf h}}.
\end{align}
Using Eq.~\eqref{eq:Rhat}, the composition of torsor functors acting on the braiding operators is given by
\begin{align}
\mathcal{F}_{\mathfrak{t}'', \Xdot''} \circ \mathcal{F}_{\mathfrak{t}', \Xdot'} \left( R^{a_{\bf g} b_{\bf h} } \right) &= \SymmAcGL
\notag \\
&=\mathscr{Q}({\bf g,h})F^{\mathfrak{q}''({\bf g,h} ) \mathfrak{q}'({\bf g,h}) a_{\bf g}'  } \SymmAcGR
\notag \\
&= \mathcal{G}_{\Gamma, \gamma} \circ \mathcal{F}_{\mathfrak{t}, \Xdot} \left( R^{a_{\bf g} b_{\bf h}} \right)
,
\label{eq:SymAcGaugT}
\end{align}
where $a_{\bf g}' = \overline{\mathfrak{q}'}({\bf g,h}) \otimes \overline{\mathfrak{q}''}({\bf g,h}) \otimes a_{\bf g}$, and $\mathscr{Q}$ is defined through the evaluation of the diagram
\begin{align}
\QLeft=\mathscr{Q}({\bf g,h}) \QRight.
\end{align}
The explicit evaluation yields
\begin{align}
\mathscr{Q}({\bf g,h})
&=\frac{1}{R^{\mathfrak{q}''({\bf g,h}) \mathfrak{t}'({\bf g,h})}}
\frac{
F^{\mathfrak{t}''({\bf g,h}) \mathfrak{q}''({\bf g,h}) \mathfrak{t}'({\bf g,h}) }
F^{\mathfrak{t}({\bf g,h}) \mathfrak{q}''({\bf g,h}) \mathfrak{q}'({\bf g,h}) }
}
{
F^{\mathfrak{t}''({\bf h, \bar{h}gh}) \mathfrak{t}'({\bf g,h})\mathfrak{q}'({\bf g,h}) }
F^{\mathfrak{t}''({\bf g,h}) \mathfrak{t}'({\bf g,h}) \mathfrak{q}''({\bf g,h})}
}
.
\end{align}
Thus, we see the coefficient of the last diagram in Eq.~\eqref{eq:SymAcGaugT} yields the symmetry action gauge transformation
\begin{align}
\gamma_{a_{\bf g}}({\bf h} ) = \mathscr{Q}({\bf g,h})F^{\mathfrak{q}''({\bf g,h} ) \mathfrak{q}'({\bf g,h}) a_{\bf g}'  }
,
\end{align}
which gives Eq.~\eqref{eq:gamma_comp}.

We note that $\mathscr{Y}_{\mathfrak{t}'', \mathfrak{t}'}$ is not necessarily equal to $\mathscr{Y}_{\mathfrak{t}', \mathfrak{t}''}$.
In particular, we have
\begin{align}
[\mathscr{Y}_{\mathfrak{t}'', \mathfrak{t}'} ] &= [\mathscr{C}_{\mathfrak{t}', \mathfrak{t}''} \mathscr{Y}_{\mathfrak{t}', \mathfrak{t}''} ],
\\
\label{eq:commutation}
\mathscr{C}_{\mathfrak{t}', \mathfrak{t}''}({\bf g,h,k}) &= \frac{M_{\mathfrak{t}'({\bf g,h}) \mathfrak{t}''({\bf gh,k}) } }
{M_{{}^{\bf g} \mathfrak{t}'({\bf h,k}) \mathfrak{t}''({\bf g,hk}) }}
,
\end{align}
where the square brackets here indicate the $H^3(G,\text{U}(1))$ cohomology classes.
We obtain this by result by rearranging the order of the torsor functors' cocyclators and braiding the necessary phantom lines to obtain the different order of torsor functor actions, giving the diagrammatic relation
\begin{align}
\CommR = \mathscr{C}_{\mathfrak{t}', \mathfrak{t}''}({\bf g,h,k}) \text{d} \Xi({\bf g,h,k})\; \CommL.
,
\end{align}
with $\Xi({\bf g,h}) = M_{\mathfrak{t}'({\bf g,h})  \mathfrak{t}''({\bf g,h}) }$.

\section{Example: Trivial Symmetry Action on \texorpdfstring{$\mathcal{C}_{\bf 0}$}{C0}}
\label{sec:TrivSymmAc}

When the symmetry action on $\mathcal{C}_{\bf 0}$ is trivial, i.e. $\rho_{\bf g} = \openone$ on $\mathcal{C}_{\bf 0}$ for all ${\bf g}$, we have a particularly simple application of the torsor functor.
In this case, there is always a (trivial) $G$-crossed extension given by $\mathcal{C}_{G}^{\times} = \mathcal{C}_{\bf 0} \boxtimes \spt_{G}^{[1]}$, which we will use as the base theory to which the torsor functor is applied.
Recall that the simple objects of $\spt_{G}^{[1]}$ are the group elements of $G$, with fusion given by group multiplication, and all the basic data are equal to 1.
Thus, the topological charges of $\mathcal{C}_{G}^{\times}$ take the form $a_{\bf g} = (a , {\bf g})$, and the basic data are independent of the $G$ labels, i.e. are given by the basic data of $\mathcal{C}_{\bf 0}$.
In this section, topological charges written without group labels are understood to be in the $\mathcal{C}_{\bf 0}$ sector.
Applying a torsor functor to $\mathcal{C}_{G}^{\times}$, we can write the basic data of the new theory $\widehat{\mathcal{C}}_{G}^{\times}=\mathcal{F}_{\mathfrak{t}, \Xdot}(\mathcal{C}_{G}^{\times})$ in terms of the basic data of $\mathcal{C}_{\bf 0}$.
In the case where $\mathcal{C}_{\bf 0}$ is a MTC, applying torsor functors to $\mathcal{C}_{G}^{\times}$ generates all possible $G$-crossed extensions of $\mathcal{C}_{\bf 0}$ with trivial symmetry action on $\mathcal{C}_{\bf 0}$, as follows from the $H^2(G,\mathcal{A})$ and $H^3(G,\text{U}(1))$ classification.
(The case where $\mathcal{C}_{\bf 0}$ is a MTC with no fusion multiplicities was obtained in Ref.~\onlinecite{Bark2019} by directly solving the $G$-crossed consistency equations; the basic data found there are presumably gauge equivalent to the data found here for those cases.)

In detail, we obtain the fusion rules
\begin{align}
a_{\bf g} \widehat{\otimes} b_{\bf h} &\cong \bigoplus_{c_{\bf gh}} N_{a b}^{\bar{\mathfrak{t}}({\bf g}, {\bf h}) \otimes c} c_{\bf gh}
, \\
\widehat{N}_{a_{\bf g} b_{\bf h}}^{c_{\bf gh}} &= N_{a b}^{\bar{\mathfrak{t}}({\bf g}, {\bf h}) \otimes c}
.
\end{align}
The symmetry action on defect charges is given by
\begin{equation}
\widehat{\rho}_{\bf k} (a_{\bf g}) =   \bar{\mathfrak{q}}({\bf g},{\bf \bar{k}} ) \otimes a_{\bf g}
.
\end{equation}

The remaining basic data is given by
\begin{align}
& \left[\widehat{F}^{a_{\bf g} b_{\bf h} c_{\bf k}} _{d_{\bf ghk}}\right]_{(e_{\bf gh},\alpha,\beta)(f_{\bf hk},\mu,\nu)} =
\sum_{\beta',\nu', \nu''}
\;
\left[{F}^{\mathfrak{t}({\bf g},{\bf h}) [\bar{\mathfrak{t}}({\bf g},{\bf h}) e] c}_{[\bar{\mathfrak{t}}({\bf gh},{\bf k}) d]}\right]_{(e,\beta), ([\bar{\mathfrak{t}}({\bf g},{\bf h}) \bar{\mathfrak{t}}({\bf gh},{\bf k}) d],\beta')}
\notag \\
& \qquad\qquad \times
\left[{F}^{a b c} _{[\bar{\mathfrak{t}}({\bf g},{\bf h}) \bar{\mathfrak{t}}({\bf gh},{\bf k})d]}\right]_{([\bar{\mathfrak{t}}({\bf g},{\bf h}) e],\alpha, \beta'), ([\bar{\mathfrak{t}}({\bf h},{\bf k}) f],\mu ,\nu'')}
\left( F^{\mathfrak{t}({\bf gh},{\bf k} ) \mathfrak{t}({\bf g},{\bf h}) [\bar{\mathfrak{t}}({\bf g},{\bf h}) \bar{\mathfrak{t}}({\bf gh},{\bf k}) d]}_{d} \right)^{-1}
\notag \\
& \qquad \qquad
\times
F^{\mathfrak{t}({\bf g},{\bf hk}) \mathfrak{t}({\bf h},{\bf k} ) [\bar{\mathfrak{t}}({\bf g},{\bf h}) \bar{\mathfrak{t}}({\bf gh},{\bf k}) d]}_{d}
\left[\left( F^{\mathfrak{t}({\bf h},{\bf k}) a [\bar{\mathfrak{t}}({\bf h},{\bf k}) f]}_{[\bar{\mathfrak{t}}({\bf g},{\bf hk})  d] }  \right)^{-1}\right]_{([\bar{\mathfrak{t}}({\bf g},{\bf h}) \bar{\mathfrak{t}}({\bf gh},{\bf k})d],\nu''), ([\mathfrak{t}({\bf h},{\bf k}) a],\nu')}
\notag \\
&
\qquad \qquad \times
\left( R^{\mathfrak{t}({\bf h},{\bf k}) a}_{[{\mathfrak{t}}({\bf h},{\bf k}) a] }\right)^{-1}
\left[ F^{a \mathfrak{t}({\bf h},{\bf k}) [\bar{\mathfrak{t}}({\bf h},{\bf k}) f]}_{[\bar{\mathfrak{t}}({\bf g},{\bf hk}) d]} \right]_{([{\mathfrak{t}}({\bf h},{\bf k}) a],\nu'), (f,\nu)}
\Xdot({\bf g},{\bf h}, {\bf k})
,
\end{align}
\begin{align}
\left[\widehat{R}^{a_{\bf g} b_{\bf h}}_{c_{\bf gh}} \right]_{\mu \nu}
=\sum_{\nu'}
\left[R^{[\bar{\mathfrak{q}}({\bf g} , {\bf h}) a] b}_{[\bar{\mathfrak{t}}({\bf h},{\bf \bar{h}gh}) c]} \right]_{\mu \nu'}
F^{\mathfrak{t}({\bf g },{\bf h})\mathfrak{q}({\bf g},{\bf h}) [\bar{\mathfrak{t}}({\bf h},{\bf \bar{h}gh}) c]}_{c}
\left[ \left(F^{\mathfrak{q}({\bf g},{\bf h}) [\bar{\mathfrak{q}}({\bf g} , {\bf h}) a] b}_{[\bar{\mathfrak{t}}({\bf g},{\bf h}) c]}\right)^{-1} \right]_{([\bar{\mathfrak{t}}({\bf h},{\bf \bar{h}gh})c], \nu')(a,\nu) }
,
\end{align}
\begin{align}
& \left[\widehat{U}_{\bf k}\left( a_{\bf g}, b_{\bf h}; c_{\bf gh} \right)\right]_{\mu \nu} =
\sum_{\nu',\nu'',\nu'''}
\left( R^{\mathfrak{q}({\bf h},{\bf k}) a} \right)^{-1}
\left[ \left(F^{\mathfrak{q}({\bf g},{\bf k}) [\bar{\mathfrak{q}}({\bf g},{\bf k}) a] [\bar{\mathfrak{q}}({\bf h},{\bf k}) b]}_{[\bar{\mathfrak{t}}({\bf g},{\bf h}) \bar{\mathfrak{q}}({\bf h},{\bf k}) c]}\right)^{-1} \right]_{([\bar{\mathfrak{t}}({\bf g},{\bf h}) \bar{\mathfrak{q}}({\bf g},{\bf k}) \bar{\mathfrak{q}}({\bf h},{\bf k}) c], \nu')(a,\nu'') }
\notag \\
& \qquad \times
\left[ \left(F^{\mathfrak{q}({\bf h},{\bf k}) a [\bar{\mathfrak{q}}({\bf h},{\bf k}) b]}_{[\bar{\mathfrak{t}}({\bf g},{\bf h}) c]}\right)^{-1} \right]_{([\bar{\mathfrak{t}}({\bf g},{\bf h}) \bar{\mathfrak{q}}({\bf h},{\bf k}) c, \nu'')([\mathfrak{q}({\bf h},{\bf k})a],\nu''') }
\left[ F^{a \mathfrak{q}({\bf h},{\bf k}) [\bar{\mathfrak{q}}({\bf h},{\bf k}) b]}_{[\bar{\mathfrak{t}}({\bf g},{\bf h}) c]} \right]_{([\mathfrak{q}({\bf h},{\bf k})a],\nu''')(b_{\bf h}, \nu) }
\notag \\
& \qquad \times
\frac{F^{[\mathfrak{t}({\bf gk},{\bf \bar{k}hk})  \mathfrak{t}({\bf g},{\bf k})] \mathfrak{q}({\bf g},{\bf k})  [\bar{\mathfrak{t}}({\bf g},{\bf h}) \bar{\mathfrak{q}}({\bf g},{\bf k}) \bar{\mathfrak{q}}({\bf h},{\bf k}) c]  }
F^{[\mathfrak{t}({\bf h},{\bf k}) \mathfrak{t}({\bf g},{\bf hk}) ] \mathfrak{q}({\bf h},{\bf k})  [\bar{\mathfrak{t}}({\bf g},{\bf h}) \bar{\mathfrak{q}}({\bf h},{\bf k}) c]}
F^{\mathfrak{t}({\bf gh},{\bf k}) \mathfrak{t}({\bf g},{\bf h})  [\bar{\mathfrak{t}}({\bf g},{\bf h}) c] }
}
{F^{\mathfrak{t}({\bf k},{\bf \bar{k}ghk}) \mathfrak{t}({\bf \bar{k}gk},{\bf \bar{k}hk})  [\bar{\mathfrak{t}}({\bf g},{\bf h}) \bar{\mathfrak{q}}({\bf g},{\bf k}) \bar{\mathfrak{q}}({\bf h},{\bf k}) c]  }
F^{\mathfrak{t}({\bf gh},{\bf k}) \mathfrak{q}({\bf gh},{\bf k})  [\bar{\mathfrak{q}}({\bf gh},{\bf k}) c]  }
}
\notag \\
& \qquad \times
\frac{
F^{\mathfrak{t}({\bf g},{\bf hk}) \mathfrak{t}({\bf h},{\bf k}) \mathfrak{q}({\bf h},{\bf k})  }
}
{
F^{\mathfrak{t}({\bf gk},{\bf \bar{k}hk}) \mathfrak{t}({\bf g},{\bf k})  \mathfrak{q}({\bf g},{\bf k})  }
}
\frac{\Xdot({\bf g},{\bf k}, {\bf \bar{k}hk}) }{\Xdot({\bf g},{\bf h}, {\bf k}) \Xdot({\bf k}, {\bf \bar{k}gk}, {\bf \bar{k}hk})}
,
\end{align}
\begin{align}
& \widehat{\eta}_{x_{\bf k}} \left( {\bf g}, {\bf h} \right) =
R^{\mathfrak{t}({\bf g},{\bf h}) x}
R^{ [\bar{\mathfrak{q}}({\bf k},{\bf gh}) x] \mathfrak{t}({\bf g},{\bf h})}
\notag \\
& \quad \times
\frac{F^{\mathfrak{t}({\bf kg},{\bf h}) \mathfrak{t}({\bf k},{\bf g}) \mathfrak{q}({\bf k},{\bf g})  }
F^{\mathfrak{t}({\bf g},{\bf \bar{g}kgh}) \mathfrak{t}({\bf \bar{g}kg},{\bf h}) \mathfrak{q}({\bf \bar{g}kg},{\bf h})  }
F^{\mathfrak{t}({\bf gh},{\bf \bar{h}\bar{g}kgh}) \mathfrak{t}({\bf g},{\bf h}) [\bar{\mathfrak{q}}({\bf k},{\bf gh}) x] }
F^{\mathfrak{t}({\bf k},{\bf gh}) \mathfrak{q}({\bf k},{\bf gh})  [{\mathfrak{t}}({\bf g},{\bf h}) \bar{\mathfrak{q}}({\bf k},{\bf gh}) x] }
}
{F^{\mathfrak{t}({\bf k},{\bf gh}) \mathfrak{t}({\bf g},{\bf h}) x }
F^{[\mathfrak{t}({\bf kg},{\bf h})  \mathfrak{t}({\bf k},{\bf g})] \mathfrak{q}({\bf k},{\bf g})  [\bar{\mathfrak{q}}({\bf k},{\bf g}) x]  }
F^{[\mathfrak{t}({\bf g},{\bf \bar{g}kgh}) \mathfrak{t}({\bf \bar{g}kg},{\bf h}) ] \mathfrak{q}({\bf \bar{g}kg},{\bf h})  [\bar{\mathfrak{q}}({\bf k},{\bf gh}) x ]   }
F^{\mathfrak{q}({\bf k},{\bf gh}) [\bar{\mathfrak{q}}({\bf k},{\bf gh}) x] \mathfrak{t}({\bf g},{\bf h}) }
}
\notag \\
& \quad \times
\frac{\Xdot({\bf g},{\bf \bar{g}kg}, {\bf h}) }{\Xdot({\bf g},{\bf h}, {\bf \bar{h}\bar{g}kgh}) \Xdot({\bf k}, {\bf g}, {\bf h})}
.
\end{align}

Finally, we note that the obstructions for $\mathcal{C}_{G}^{\times}$ are trivial, so the defectification obstruction for $\widehat{\mathcal{C}}_{G}^{\times}$ is given by the relative obstruction
\begin{align}
\widehat{\defectO}({\bf g},{\bf h},{\bf k},{\bf l}) &= \defectO_{r}(\mathfrak{t})({\bf g},{\bf h},{\bf k},{\bf l}) \notag \\
&= R^{\mathfrak{t}({\bf k},{\bf l}) \mathfrak{t}({\bf g},{\bf h})}
\frac{
F^{\mathfrak{t}({\bf gh},{\bf kl}) \mathfrak{t}({\bf g},{\bf h}) \mathfrak{t}({\bf k},{\bf l})}
F^{\mathfrak{t}({\bf g},{\bf hkl}) \mathfrak{t}({\bf hk},{\bf l} ) \mathfrak{t}({\bf h},{\bf k})}
F^{\mathfrak{t}({\bf ghk},{\bf l}) \mathfrak{t}({\bf gh},{\bf k} ) \mathfrak{t}({\bf g},{\bf h})}
}
{
F^{\mathfrak{t}({\bf gh},{\bf kl}) \mathfrak{t}({\bf k},{\bf l})\mathfrak{t}({\bf g},{\bf h})}
F^{\mathfrak{t}({\bf g},{\bf hkl}) \mathfrak{t}({\bf h},{\bf kl}) \mathfrak{t}({\bf k},{\bf l}) }
F^{\mathfrak{t}({\bf ghk},{\bf l}) \mathfrak{t}({\bf g},{\bf hk}) \mathfrak{t}({\bf h},{\bf k} ) }
}
.
\end{align}

\begin{acknowledgments}
We are grateful to Corey Jones, Ryan Thorngren, and Zhenghan Wang for helpful conversations.
C.K. was supported in part by the Institute for Quantum Information and Matter, an NSF Frontier center funded by the Gordon and Betty Moore Foundation, as well as by the Walter Burke Institute for Theoretical Physics at Caltech.
\end{acknowledgments}

\appendix

\section{Background and notation}
\label{sec:background}

This appendix reviews braided tensor categories (BTCs) and their $G$-crossed extensions.
Throughout, we employ a diagrammatic formalism to encapsulate complicated calculations in a simple way.
Additional details can be found in Ref.~\onlinecite{Bark2019}.

\subsection{Braided tensor categories}
\label{sec:MTC}

Let $\bMTC$ be a BTC with simple objects (topological charges) $a,b,\ldots \in \bMTC$.
The simple objects satisfy an associative fusion algebra
\begin{align}
\label{eq:fusion}
a\otimes b &= \bigoplus_{c\in \bMTC} N_{ab}^c c
\end{align}
where the fusion multiplicities $N_{ab}^c$ are non-negative integers indicating the number of distinct ways that $a$ and $b$ can be combined to form $c$.

There is a unique topological vacuum charge, denoted $\I $ for all theories.
In the diagrammatic formalism, we can freely add and remove charge lines associated with $\I $.
Each topological charge $a$ has a unique conjugate charge denoted $\bar{a}$ such that the fusion of $a$ and $\bar{a}$ includes $\I $ once, i.e. $N_{a b}^{\I } = \delta_{b \bar{a}}$.
It is possible for charges to be self-dual, meaning $a=\bar{a}$.
When $\sum_c N_{a \bar{a}}^c=1$, then $a$ is {\it Abelian}.
If instead $\sum_c N^c_{a \bar{a}}\neq 1$, then $a$ is {\it non-Abelian}.
The vacuum charge $\I $ is necessarily Abelian, $\I \otimes a = a$, thus all MTCs $\bMTC$ contain at least one (trivial) Abelian charge.
We define the {\it quantum dimension} $d_a$ to be the largest eigenvalue of the fusion matrix $N_{a}$ defined as $\left( N_a\right)_b^c = N_{ab}^c$, with $b$ and $c$ the matrix indices.
Abelian charges always have $d_a=1$, while non-Abelian charges have $d_a>1$.
The quantum dimensions satisfy $d_a d_b = \sum_{c\in \bMTC} N_{ab}^c d_c$, from which we clearly see that two Abelian charges always fuse to an Abelian charge.
Together with each charge having a unique conjugate charge, this implies that the collection of Abelian charges $\mathcal{A}$ of any MTC $\bMTC$ forms a group under fusion.
The total quantum dimension $\mathcal{D}$ of the MTC $\bMTC$ is given by $\mathcal{D}=\sqrt{\sum_{a\in \bMTC} d_a^2}$.

The fusion of $a$ and $b$ to $c$ has an associated vector space $V_{ab}^c$ with $\text{dim}V_{ab}^c= N_{ab}^c$, and a dual (splitting) space $V_c^{ab}$.
In the diagrammatic notation, topological charges label oriented line segments and trivalent vertices describe states within these vector spaces
\begin{align}
\VabcmuR &\propto \bra{a,b;c,\mu} \in V_{ab}^c,&
\VabcmuL &\propto \ket{a,b;c,\mu} \in V_c^{ab}
\end{align}
where $\mu =1,\dots,N_{ab}^c$.
Here we will focus on theories without multiplicity, so that $N_{ab}^c =0,1$ and the fusion/splitting spaces are one-dimensional.
We use the conventions common in the physics literature where,
\begin{align}
\idab = \sum_{c,\mu} \sqrt{\frac{d_c}{d_a d_b}}\; \idabesolve
\end{align}
and the bubble removal identity,
\begin{align}
\bubbleabc = \delta_{cc'}\delta_{\mu \mu'}\sqrt{\frac{d_a d_b}{d_c}}\;\idc.
\end{align}

\begin{figure}[t!]
   \centering
   \includegraphics[width=.6\columnwidth]{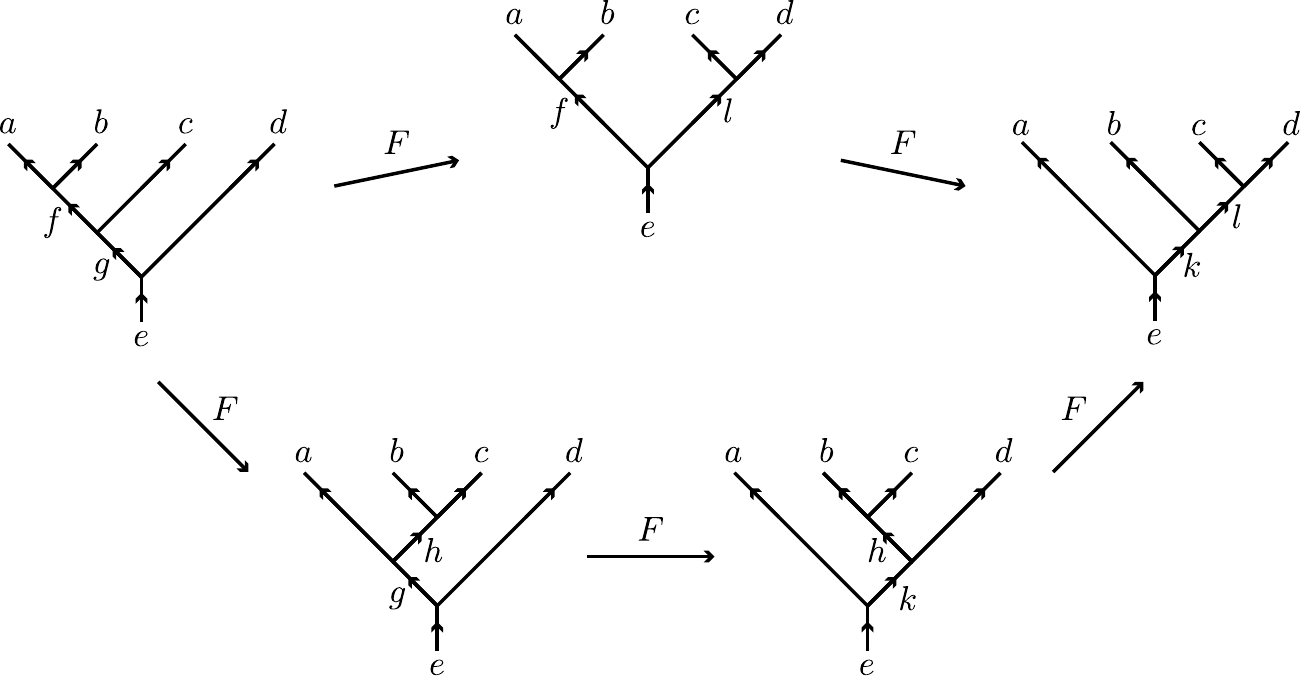}
   \caption{The pentagon equation is a consistency condition that requires equivalence of different sequences of $F$-moves that start and end at common decompositions of the fusion state space.}
   \label{fig:pentagon}
\end{figure}

More complicated diagrams involving additional charge lines are constructed by stacking trivalent vertices and connecting lines corresponding to the same topological charge.
The resulting states belong to the fusion/splitting spaces of three or more topological charges.
These vector spaces satisfy associativity of fusion
\begin{align}
V_d^{abc} &\cong \bigoplus_e V_e^{ab} \otimes V_d^{ec} \cong \bigoplus_f V_d^{af} \otimes V_f^{bc},
\end{align}
where $\cong$ denotes an isomorphism called an $F$-move, written diagrammatically as
\begin{align}
\FLeftMTC &= \sum_{f,\mu, \nu} [F_d^{abc}]_{(e,\alpha, \beta)(f,\mu, \nu)} \FRightMTC.
\end{align}
The $F$-moves amount to changing the basis of the splitting space of three topological charges; as such they correspond to unitary transformations.
Unitarity fixes
\begin{align}
\left[ (F_d^{abc})^{-1}\right]_{(f,\mu,\nu)(e,\alpha,\beta)} &= \left[(F_d^{abc})^\dagger \right]_{(f,\mu,\nu)(e,\alpha,\beta)} = \left[ F_d^{abc}\right]^*_{(e,\alpha,\beta)(f,\mu, \nu)}.
\end{align}
Any combination of $F$-moves that begin and end with the same diagram must be equivalent; this consistency condition results in the pentagon equation for the $F$-symbols:
\begin{align}
\label{eq:pentagoneqn}
&\sum_\delta [F_e^{fcd}]_{(g,\beta,\gamma),(l,\nu,\delta)} [F_e^{abl}]_{(f,\alpha,\delta),(k,\mu,\lambda)} \\
\notag
&\quad \quad \quad = \sum_{h,\sigma,\psi,\rho} [F_g^{abc}]_{(f,\alpha,\beta),(h,\psi,\sigma)} [F_e^{ahd}]_{(g,\sigma,\gamma),(k,\rho,\lambda)} [F_k^{bcd}]_{(h,\psi,\rho),(l,\nu,\mu)}.
\end{align}
The corresponding diagrammatic equation is depicted in Fig.~\ref{fig:pentagon}.

Braiding refers to the exchange of a pair of topological charges.
At the level of topological charges, the interplay of locality with braiding requires the fusion rules to be commutative, i.e. $a \otimes b = b \otimes a$.
We describe the action of braiding on the fusion and splitting state space using $R$-moves that map between vector spaces as $R^{ab}:V_c^{ba}\to V_c^{ab}$.
A counterclockwise braid is represented diagrammatically as
\begin{align}
R^{ab} = \Rab
\end{align}
It is convenient to write out the matrix elements of $R^{ab}$ diagrammatically as,
\begin{align}
\Rabc =\sum_\nu [R^{ab}_c]_{\mu \nu }\RVabc.
\end{align}
An important related object is the double braid, or {\it monodromy}, of two charges $a$ and $b$.
If at least one of $a$ or $b$ is Abelian, then the monodromy $M_{ab}$ is defined as the phase relating the double braid of $a$ with $b$ to the identity operator
\begin{align}
R^{ab} R^{ba} &= \Mab = M_{ab}\;\idab.
\end{align}

\begin{figure*}[t!]
   \centering
   \includegraphics[width=.99\columnwidth]{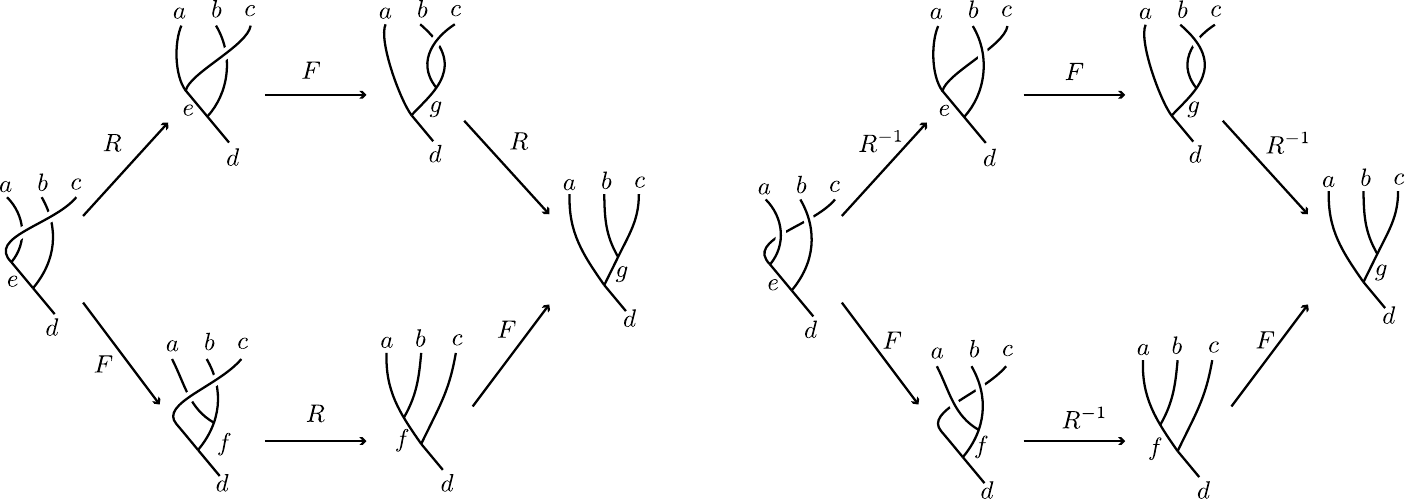}
   \caption{Consistency of braiding and fusion is described by the hexagon equation.
Implicit in the braiding consistency is the property that a line can freely slide over or under a vertex, which follows from locality.
   }
   \label{fig:hexagon}
\end{figure*}

Compatibility of braiding with fusion implies that a line may freely slide over or under a vertex; this results in a pair of consistency conditions known as the hexagon equations
\begin{align}
&\sum_{\lambda,\gamma} [R_e^{ac}]_{\alpha,\lambda}\, [F_d^{acb}]_{(e,\lambda,\beta),(g,\gamma,\nu)} \,[R_g^{bc}]_{\gamma \mu}\\ \notag
&\quad = \sum_{f,\sigma,\delta,\psi} [F_d^{cab}]_{(e,\alpha,\beta),(f,\delta,\sigma)} \,[R_d^{fc}]_{\sigma\psi} \, [F_d^{abc}]_{(f,\delta,\psi),(g,\mu,\nu)}
\end{align}
and
\begin{align}
&\sum_{\lambda,\gamma}\left[(R_e^{ca})^{-1}\right]_{\alpha \lambda} [F_d^{acb}]_{(e,\lambda,\beta),(g,\gamma,\nu)} \left[(R^{cb}_g)^{-1}\right]_{\gamma \mu} \\
&\quad \notag =\sum_{f,\sigma,\delta,\psi} [F_d^{cab}]_{(e,\alpha,\beta),(f,\delta,\sigma)}\left[ (R_d^{cf})^{-1}\right]_{\sigma \psi} [F_d^{abc}]_{(f,\delta,\psi),(g,\mu,\nu)}
\end{align}
depicted diagrammatically in Fig.~\ref{fig:hexagon}.
The $F$- and $R$-symbols constitute the basic data characterizing a BTC.

Unitary transformations $\Gamma_c^{ab}$ acting on the fusion/splitting spaces $V_c^{ab}$ and $V_{ab}^c$ redefine the basis states as
\begin{align}\label{eq:vbgt}
\widetilde{\ket{a,b;c,\mu}} &= \sum_\nu \left[\Gamma_c^{ab}\right]_{\mu \nu} \ket{a,b;c,\nu}.
\end{align}
For the case of multiplicity free theories, these are complex phases that define a gauge freedom of the theory.
The $F$- and $R$-symbols are not invariant under such {\it vertex basis gauge transformations}, but rather transform as
\begin{align}
\left[\widetilde{F}_d^{abc} \right]_{(e,\alpha,\beta)(f,\mu,\nu)} &= \sum_{\alpha',\beta',\mu',\nu'}
[\Gamma_e^{ab}]_{\alpha \alpha'}
[\Gamma_c^{ec}]_{\beta \beta'}
 \left[F^{abc}_d \right]_{(e,\alpha',\beta')(f,\mu',\nu')}
 [( \Gamma_f^{bc})^{-1} ]_{\mu' \mu}
[( \Gamma_d^{af})^{-1} ]_{\nu' \nu}
\\
\left[ \widetilde{R}_c^{ab} \right]_{\mu \nu}&= \sum_{\mu',\nu'}\left[\Gamma_c^{ba}\right]_{\mu \mu'} \left[R_c^{ab} \right]_{\mu'\nu'}\left[\left(\Gamma^{ab}_c\right)^{-1}\right]_{\nu' \nu}.
\end{align}
Preserving the condition that vacuum lines may be freely added and removed fixes ${\Gamma_a^{a\I } = \Gamma_b^{\I b} = \Gamma_{\I }^{\I \I }}$.

\subsection{\texorpdfstring{$G$}{G}-crossed braided tensor categories}
\label{sec:G-BTC}

We review the $G$-crossed formalism following Ref.~\onlinecite{Bark2019} for a BTC $\bMTC$ with global symmetry $G$.
Each step in the $G$-crossed formalism fixes a new attribute of the theory (symmetry action $\to$ symmetry fractionalization $\to$ full defect theory) that must be made consistent with the previously specified attribute; if the corresponding consistency conditions cannot be satisfied, the theory is obstructed.
When these obstructions vanish, $\bMTC$ can be extended to a $G$-crossed BTC $\bMTC_G^\times$.
Classifying $G$-crossed BTCs effectively amounts to solving the $G$-crossed consistency conditions, which quickly become intractabel for general $\bMTC$ and $G$.
This process is substantially simplified by connecting the $G$-crossed consistency conditions to cohomology groups~\cite{Etingof2010}, see Appendix~\ref{app:cohomology} for a review.
In the physics community, $G$-crossed extensions generalize the algebraic theory of anyons to include symmetry defects, thereby describing $(2+1)$D symmetry enriched and symmetry protected topological phases.

A {\it topological symmetry} of a BTC $\bMTC$ is an auto-equivalence map ${\varphi:\bMTC\to \bMTC}$ that leaves the topological data invariant.
Maps $\varphi$ and $\varphi'$ related by a natural isomorphism (a vertex basis gauge transformation of the form $\Gamma^{ab}_c = \gamma_a \gamma_b /\gamma_c$) belong to the same equivalence class  $[\varphi]$.
The equivalence classes form a group $\text{Aut}(\bMTC)$.
Prominent physics examples of topological symmetries include the electric-magnetic duality in the toric code and layer-exchange symmetries in $(mml)$ quantum Hall systems.

The global symmetry $G$ is incorporated through a group action (homomorphism) ${[\rho]:G\to \text{Aut}(\bMTC)}$, which maps each group element ${\bf g}\in G$ to a topological symmetry $[\rho_{\bf g}]$ such that: $[\rho_{\bf g}]\circ [\rho_{\bf h}] = [\rho_{\bf gh}]$.
We denote the symmetry action on a topological charge $a\in \bMTC$ using the shorthand notation
\begin{align}\label{eq:shorthand}
{}^{\bf g}a &= \rho_{\bf g}(a), & {\bf \bar{g} }&= {\bf g}^{-1}.
\end{align}
Consistency of fusion requires,
\begin{align}
N_{ab}^c = \rho_{\bf g}(N_{ab}^c) = N_{{}^{\bf g}a {}^{\bf g}b}^{{}^{\bf g}c}.
\end{align}
The particular element $\rho_{\bf g} \in [\rho_{\bf g}]$ does not change the action on the topological charge or fusion rules.

\begin{figure}[t!]
\begin{center}
\includegraphics[width=.8\columnwidth]{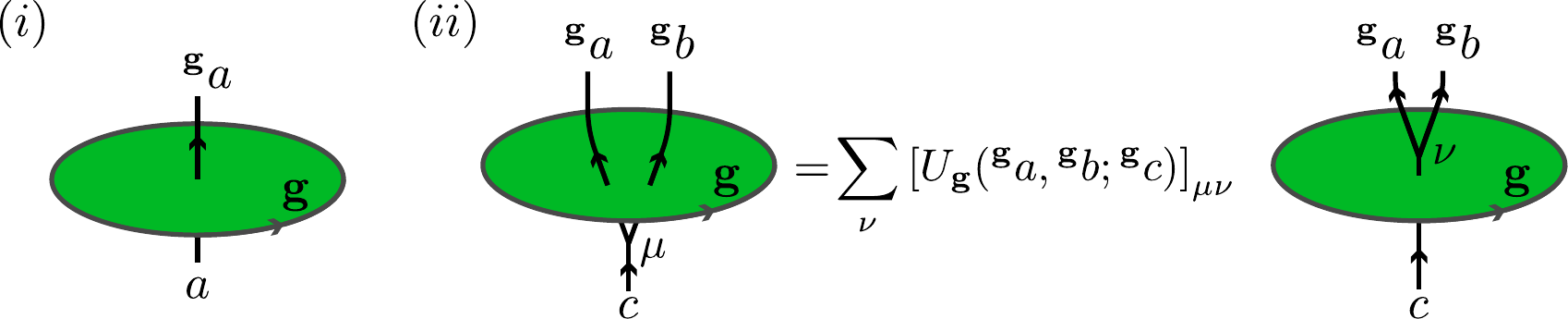}
\caption{
Green sheets represent ${\bf g}$-labeled branch sheets.
$(i)$ A simple object passing through a ${\bf g}$-sheet is acted on by $\rho_{\bf g}$.
$(ii)$ A ${\bf g}$-sheet acts non-trivially on trivalent vertices.
}
\label{fig:symmetry-action}
\end{center}
\end{figure}

In the diagrammatic formalism, we can think of the symmetry action as corresponding to a sheet, labeled by ${\bf g}$, oriented perpendicular to the vertical direction; a topological charge line piercing such a sheet is acted on by $\rho_{\bf g}$, see Fig.~\ref{fig:symmetry-action}$(i)$.
The trivial group element ${\bf 0} \in G$ always corresponds to the trivial element of $\text{Aut}(\bMTC)$, $\rho_{\bf 0}(a)=a$, thus ${\bf 0}$-sheets can be added or removed at will.

In addition to the charge permutations, the symmetry action is specified by unitary transformations acting on the splitting and fusion spaces, which we specify as:
\begin{align}
\rho_{\bf g}(\ket{a,b;c,\mu}) &=\sum_\nu \left[ U_{\bf g}({}^{\bf g}a, {}^{\bf g}b; {}^{\bf g}c)\right]_{\mu \nu} \ket{{}^{\bf g} a, {}^{\bf g}b; {}^{\bf g}c,\nu},
\end{align}
with diagrammatic representation shown in Fig.~\ref{fig:symmetry-action}$(ii)$.
The $U$-symbols correspond to the vertex basis gauge transformations needed to make the basic data ($F$-symbols and $R$-symbols) match the original data exactly after the topological charges are permuted by the symmetry action.
Consequently, for $\rho$ to be an automorphism of $\bMTC$ we require,
\begin{align}
[F^{abc}_d]_{(e,\alpha,\beta),(f,\mu,\nu)} &= \rho_{\bf g}\left([F^{abc}_d]_{(e,\alpha,\beta),(f,\mu,\nu)} \right) \notag
\\ &=\sum_{\alpha', \beta',\mu',\nu'}{[U_{\bf g}({}^{\bf g}a, {}^{\bf g} b; {}^{\bf g} e)]_{\alpha\alpha'}[U_{\bf g}({}^{\bf g}e, {}^{\bf g} c; {}^{\bf g}d)]_{\beta\beta'}}
[F^{{}^{\bf g} a{}^{\bf g} b{}^{\bf g} c}_{{}^{\bf g} d}]_{({}^{\bf g} e,\alpha',\beta'),( {}^{\bf g} f,\mu',\nu')} \notag \\
&\quad \quad \quad \times{[U_{\bf g}({}^{\bf g} a, {}^{\bf g} f; {}^{\bf g}d)^{-1}]_{\nu'\nu}[ U_{\bf g}({}^{\bf g}b, {}^{\bf g}c; {}^{\bf g}f)^{-1}]_{\mu' \mu}}
\\
[R^{ab}_c]_{\mu \nu} &= \rho_{\bf g}\left([R^{ab}_c ]_{\mu \nu}\right) = {[U_{\bf g}({}^{\bf g}b,{}^{\bf g}a; {}^{\bf g}c)]_{\mu\mu'}}\left[ R^{{}^{\bf g}a{}^{\bf g}b}_{{}^{\bf g}c}\right]_{\mu'\nu'}
{[U_{\bf g}({}^{\bf g}a,{}^{\bf g}b;{}^{\bf g}c)^{-1}]_{\nu'\nu}}.
\end{align}
The specific choice of $U$-symbols correspond to a particular choice of $\rho_{\bf g} \in [\rho_{\bf g}]$; changing $\rho_{\bf g}$ by a natural isomorphism corresponds to a gauge freedom of the symmetry action in the theory.
We address this symmetry action gauge transformation at the end of this section.

\begin{figure}[t!]
   \centering
   \includegraphics[width=.5\columnwidth]{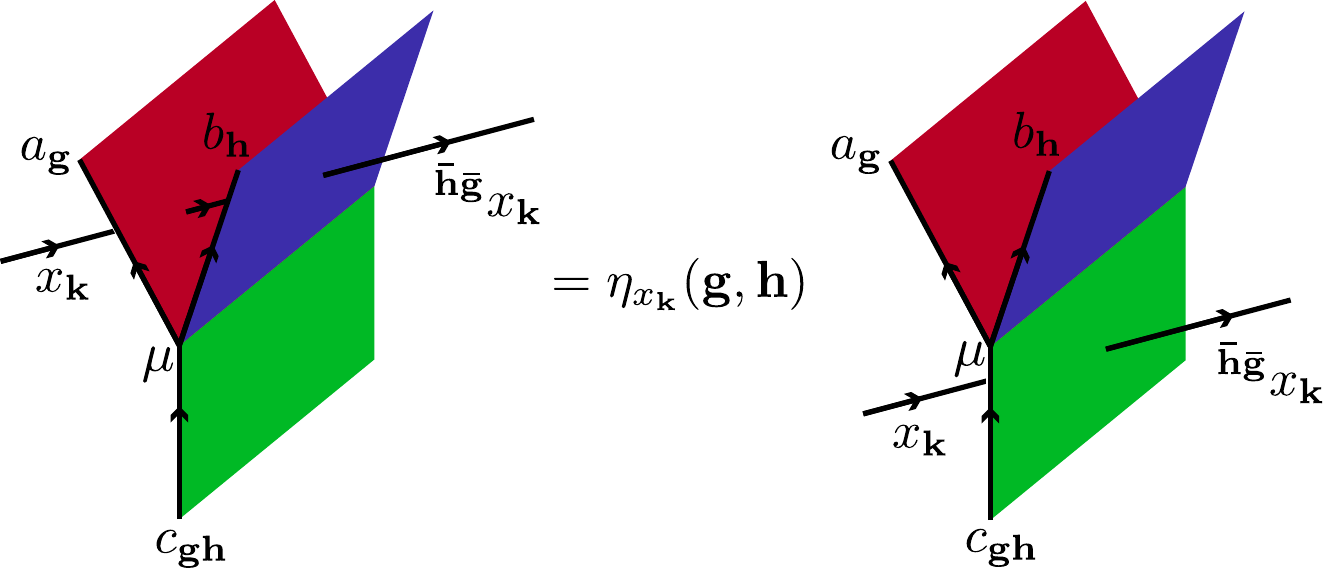}
   \caption{In this figure we show the diagrammatic interpretation of the symmetry fractionalization.
   We have suppressed the ${\bf k}$-sheets emanating off from $x_{\bf k}$ for clarity.}
   \label{fig:fractionalization}
\end{figure}

The next ingredient is the {\it symmetry fractionalization}, which specifies how topological charges carry fractionalized symmetry quantum numbers.
By assuming that the global symmetry $G$ acts on the system in an on-site manner, we can split the global symmetry action $R_{\bf g}$ into the topological symmetry action $\rho_{\bf g}$ on the anyonic state space (captured by the phases $\{U_{\bf g}(a,b;c)\}$), and a localized symmetry operation acting in the vicinity of each topological charge in the system (e.g. carried by quasiparticles)~\cite{Bark2019}.
The group multiplication structure of the latter can be reduced to a collection of projective phases $\{\eta_a({\bf g}, {\bf h})\}$ that equate the localized ${\bf gh}$ action on $a$ with successively applied ${\bf g}$ and ${\bf h}$ localized actions on $a$, see Fig.~\ref{fig:fractionalization}.

Compatibility of symmetry fractionalization with associativity of the symmetry action imposes the constraint
\begin{align}\label{eq:eta1}
{\eta}_a({\bf g}, {\bf h}) \eta_a({\bf gh},{\bf k}) &= \eta_a({\bf g},{\bf hk}) \eta_{{}^{\bar{\bf g}}a}({\bf h},{\bf k}).
\end{align}
Additionally, the $U$s and $\eta$s must satisfy the consistency condition
\begin{align}\label{eq:sliding}
\frac{\eta_a({\bf g},{\bf h}) \eta_b( {\bf g},{\bf h})  }{\eta_c({\bf g},{\bf h})}\delta_{\mu \nu} &=\sum_{\alpha,\beta}
\left[ U_{\bf g}(a,b;c)^{-1}\right]_{\mu \alpha}
\left[U_{\bf h}({}^{\bar{\bf g}}a, {}^{\bar{\bf g}}b; {}^{\bar{\bf g}}c)^{-1}\right]_{\alpha\beta}
\left[U_{\bf gh}(a,b;c)\right]_{\beta \nu}\\
&= \kappa_{\bf g,h}(a,b;c) \delta_{\mu \nu}.
\end{align}
Here, $\kappa_{\bf g,h}$ is the natural isomorphism associated with matching up the representative autoequivalence maps $\kappa_{\bf g,h} \circ \rho_{\bf g} \circ \rho_{\bf h} = \rho_{\bf gh}$.
We generalize Eqs.~\eqref{eq:eta1} and \eqref{eq:sliding} when we incorporate symmetry defects.

It is not always possible for a symmetry action $\rho$ to fractionalize consistently, meaning it may not be possible to define a set of $\eta$-symbols satsifying Eqs.~\eqref{eq:eta1} and \eqref{eq:sliding} for a fixed set of $U$-symbols.
This obstruction to fractionalization is given by an invariant $[\coho{O}]\in H^3_{[\rho]}(G,\mathcal{A})$ that indicates whether $\kappa_{\bf g,h}$ can be factored into local phases on the topological charges that are consistent with the localized symmetry action.
(For more details, see Refs.~\onlinecite{Etingof2010,Bark2019}.)
The Abelian group $\mathcal{A}$ is defined by the set of Abelian topological charges of $\bMTC$, with group multiplication specified by their fusion rules.
When $[\coho{O}]$ corresponds to the trivial element, the obstruction vanishes and the distinct fractionalization patterns form an $H^2_{[\rho]}(G,\mathcal{A})$ torsor.
In particular, the fractionalization classes do not themselves correspond to elements $[\mathfrak{t}] \in H^2_{[\rho]}(G,\mathcal{A})$, but rather are related to each other by these 2-cocycles.
A given set of $\eta$-symbols (for specified $U$-symbols) fixes the symmetry fractionalization class.
The torsorial action of $H^2_{[\rho]}(G,\mathcal{A})$ is the key focus of this paper; we return to it in Section~\ref{sec:torsor_data} and~\ref{sec:PW}.

After fixing the fractionalization class, we extend $\bMTC$ to a $G$-graded fusion category ${\bMTC_G=\bigoplus_{\bf g} \bMTC_{\bf g}}$, where the ${\bf 0}$-sector corresponds to the original theory, $\bMTC=\bMTC_{\bf 0}$, and the simple objects of $\bMTC_{\bf g}$ are symmetry defects associated with the group element ${\bf g}$.
We denote a ${\bf g}$-defect by $a_{\bf g}\in \bMTC_{\bf g}$.
This notation indicates that there may be topologically distinct objects that behave as ${\bf g}$-defects; it does {\it not } imply that a defect should generally be considered a composite object formed from binding a ${\bf g}$-flux to a topological charge $a \in \bMTC$.
The $G$-grading of the fusion category says that the fusion of defects respects the group multiplication, that is a ${\bf g}$-defect fuses with a ${\bf h}$-defect to form a ${\bf gh}$-defect:
\begin{align}
a_{\bf g}\otimes b_{\bf h} = \bigoplus_{c_{\bf gh}\in \bMTC_{\bf gh}} N_{a_{\bf g} b_{\bf h}}^{c_{\bf gh}} c_{\bf gh}
.
\end{align}

As is the case with the simple objects $a\in \bMTC$, symmetry defects must satisfy an associativity isomorphism characterized by
\begin{align}
\FLeft \!\!\!= \sum_{f_{\bf hk},\mu,\nu}
\left[ F^{a_{\bf g},b_{\bf h} c_{\bf k}}_{d_{\bf ghk} }  \right]_{(e_{\bf gh},\alpha,\beta)(f_{\bf hk},\mu,\nu)}\!\!\!\!\!
\FRight
,
\end{align}
for which the $F$-symbols must again satisfy the pentagon equation of Eq.~(\ref{eq:pentagoneqn}).
Aside from objects carrying group labels that are respected by fusion, the usual properties of fusion categories are satisfied.

The final ingredient is to extend the $G$-graded fusion category ${\bMTC_G}$ to a $G$-crossed BTC $\mathcal{C}_G^\times$, which incorporates the symmetry action into the defect braiding.
We extend the symmetry action to defects, using the shorthand notation of Eq.~\eqref{eq:shorthand}
\begin{align}
{}^{\bf h} a_{\bf g} &= \rho_{\bf h}(a_{\bf g}) \in \bMTC_{\bf hg\bar{h}}
\\  {}^{\bf h} {\bf g} &={\bf hg \bar{h}}.
\end{align}

In $(2+1)$D, a ${\bf g}$-defect worldline corresponds to a termination of the ${\bf g}$-branch sheet; in the diagrammatic formalism, we leave the ${\bf g}$-sheet implicit and assume that it emanates from the line labeled by $a_{\bf g}$ back into the page.
In this way, when a ${\bf g}$-defect line passes in front of an object, it acts on it with the ${\bf g}$-symmetry action.
This requires a $G$-crossed braiding operator
\begin{align}
R^{a_{\bf g} b_{\bf h}} = \braidingop.
\end{align}
As such, the operator $R^{a_{\bf g} b_{\bf h}}$ maps the vector space $V_{c_{\bf gh}}^{b_{\bf h} {}^{\bf \bar{ h}}a_{\bf g}}$ to $V_{c_{\bf gh}}^{a_{\bf g} b_{\bf h}}$.
We can define the corresponding $R$-symbols diagrammatically as
\begin{align}
\RLeft = \sum_\nu \left[R^{a_{\bf g} b_{\bf h}}_{c_{\bf gh}} \right]_{\mu \nu} \RRight
.
\end{align}

Additionally, in contrast to topological charge lines in $\bMTC$, defect lines cannot slide freely over vertices.
The global symmetry action phases appear as
\begin{align}
\ULeft = \sum_{\nu }\left[U_{\bf k}(a_{\bf g}, b_{\bf h}; c_{\bf gh}) \right]_{\mu \nu}\; \URight
\end{align}
while the projective local symmetry action phases arise from
\begin{align}
\etaLeft = \eta_{x_{\bf k}}({\bf g},{\bf h})\;\etaRight
\end{align}
The connection should be clear in comparison to Figs.~\ref{fig:symmetry-action} and \ref{fig:fractionalization}.
It is convenient to make the canonical gauge choices
\begin{align}
\label{Uetatrivial}
\left[ U_{\bf 0}(a_{\bf g},b_{\bf h}; c_{\bf gh}) \right]_{\mu \nu}&=\delta_{\mu \nu}\\
\eta_{\I_{\bf 0}}({\bf g},{\bf h})&=1
\end{align}
which say that sliding the trivial branch sheet over a vertex should be trivial, and similarly sliding the vacuum charge through a trivalent junction of branch sheets should also be trivial.
Requiring that vacuum lines can be freely added and removed from diagrams similarly imposes the following conditions
\begin{align}
U_{\bf k}(a_{\bf g},\I_{\bf 0}; a_{\bf g}) &= U_{\bf k}(\I_{\bf 0},b_{\bf h}; b_{\bf h}) = 1 \\\eta_{c_{\bf k}}({\bf 0},{\bf h}) &= \eta_{c_{\bf k}}({\bf g},{\bf 0}) = 1.
\end{align}

\begin{figure*}[t!]
   \centering
   \includegraphics[width=.99\textwidth]{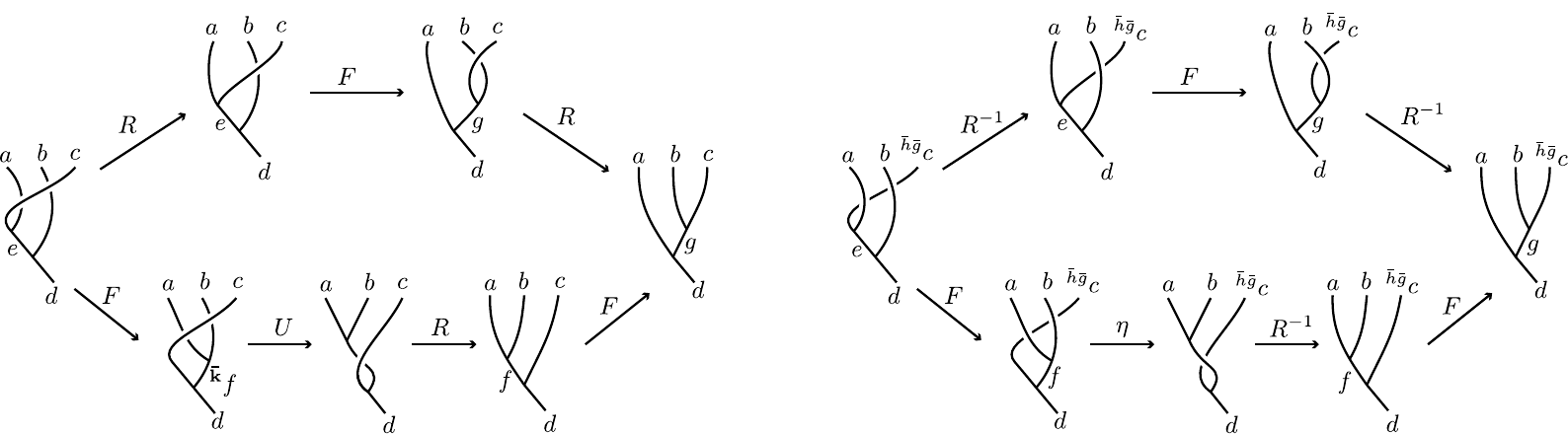} 
   \caption{Heptagon equations.
   For convenience we have left the group labels implicit, they are given by $a_{\bf g}$, $b_{\bf h}$, $c_{\bf k}$, $d_{\bf ghk}$, $e_{\bf gk}$, $f_{\bf gh}$ and $g_{\bf hk}$ for the left panel, and $a_{\bf g}$, $b_{\bf h}$, $c_{\bf k}$, $d_{\bf kgh}$, $e_{\bf kg}$, $f_{\bf gh}$, and $g_{\bf \bar{g}kgh}$ for the right panel.
   }
   \label{fig:heptagon}
\end{figure*}

Compatibility of the $R$-symbols with $F$-, $U$-, and $\eta$-symbols results in consistency conditions known as the {\it heptagon equations}.
Consistency of fusion with counterclockwise oriented braids imposes
\begin{align}
&
\sum_{\lambda,\gamma}[R_{e}^{ac}]_{\alpha \lambda}
\left[F_{d}^{ac \,^{\bf \bar{k}}b}\right] _{(e,\lambda,\beta),(g,\gamma,\nu)}
[R_{g}^{bc}]_{\gamma \mu} \notag
\\ &\quad =\sum_{f,\sigma,\delta, \eta, \psi}
\left[ F_{d}^{c \,^{\bf \bar{k}}a \,^{\bf \bar{k}}b}\right] _{( e,\alpha,\beta),(\,^{\bf \bar{k}}f,\delta,\sigma)}
[ U_{\bf k}\left( a ,b ;f \right)]_{\delta,\eta}
\left[R_{d}^{fc}\right]_{\sigma \psi}
\left[F_{d}^{abc}\right] _{(f,\eta,\psi),(g,\mu,\nu)}
,
\label{eq:heptagon+}
\end{align}
with labels $a_{\bf g}$, $b_{\bf h}$, $c_{\bf k}$, $d_{\bf ghk}$, $e_{\bf gk}$, $f_{\bf gh}$ and $g_{\bf hk}$.
Consistency of fusion with clockwise braids fixes
\begin{align}
&
\sum_{\lambda,\gamma}\left[ \left( R_{e}^{ca}\right) ^{-1}\right]_{\alpha \lambda}
\left[ F_{d}^{a \,^{\bf \bar{g}}cb}\right] _{(e,\lambda,\beta),(g,\gamma,\nu)}
\left[ \left( R_{g}^{\,^{\bf \bar{g}}c b}\right) ^{-1}\right]_{\gamma \mu} \notag
\\ &\quad =\sum_{f,\sigma,\delta,\psi}
\left[ F_{d}^{cab}\right] _{(e,\alpha,\beta),(f,\delta,\sigma)}
\eta_{c}\left({\bf g},{\bf h}\right)
\left[\left(R_{d}^{cf}\right) ^{-1}\right]_{\sigma,\psi}
\left[ F_{d}^{ab\,^{\bf \bar{h} \bar{g}}c}\right]_{(f,\delta,\psi),(g,\mu,\nu)}
\label{eq:heptagon-}
\end{align}
with labels $a_{\bf g}$, $b_{\bf h}$, $c_{\bf k}$, $d_{\bf kgh}$, $e_{\bf kg}$, $f_{\bf gh}$, and $g_{\bf \bar{g}kgh}$.
The heptagon equations are depicted diagrammatically in Fig.~\ref{fig:heptagon} with the group labels on topological charges left implicit, i.e. $a_{\bf g}$, $b_{\bf h}$, $c_{\bf k}$, etc.

Clearly the heptagon equations are generalizations of the hexagon equations that account for the nontrivial symmetry action when sliding a charge line over or under a vertex.

We now return to the gauge freedom introduced by the equivalence of autoequivalence maps differing by a natural isomorphism $\Upsilon$, which is defined as
\begin{align}
\Upsilon(a_{\bf g})&= a_{\bf g}
\\ \Upsilon(\ket{a_{\bf g},b_{\bf h} ;c_{\bf gh},\mu })&= \frac{\gamma_{a_{\bf g}} \gamma_{b_{\bf h}} }{\gamma_{c_{\bf gh}} }\ket{a_{\bf g} ,b_{\bf h} ;c_{\bf gh},\mu}
\end{align}
for some phases $\gamma_{a_{\bf g}}$.
We call the equivalence classes of autoequivalence maps ``topological symmetries.''
In this way, we consider symmetry actions $\rho$ and $\check{\rho}$ to be equivalent when
\begin{align}
\label{eq:symm_act_gauge}
\check{\rho}_{\bf g} &= \Upsilon_{\bf g} \circ \rho_{\bf g}
,
\end{align}
and we write the equivalence class as $[\rho]$.

The symmetry action gauge transformation of Eq.~\eqref{eq:symm_act_gauge} acts on the basic data according to
\begin{align}
\check{F}^{a_{\bf g} b_{\bf h} c_{\bf k}}_{d_{\bf ghk}} &= F^{a_{\bf g} b_{\bf h} c_{\bf k}}_{d_{\bf ghk}},
\\ \check{R}_{c_{\bf gh}}^{a_{\bf g} b_{\bf h}} &= \gamma_{a_{\bf g}}({\bf h}) R_{c_{\bf gh}}^{a_{\bf g} b_{\bf h}},
\\ \check{U}_{\bf k}(a_{\bf g},b_{\bf h};c_{\bf gh}) &= \frac{\gamma_{a_{\bf g}}({\bf k})\gamma_{b_{\bf h}}({\bf k})}{\gamma_{c_{\bf gh}}({\bf k})} U_{\bf k}(a_{\bf g},b_{\bf h};c_{\bf gh}),
\\ \check{\eta}_{c_{\bf k}}({\bf g},{\bf h}) &= \frac{\gamma_{c_{\bf k}}({\bf gh})}{\gamma_{{}^{\bar{\bf g}}c_{\bf k}}({\bf h}) \gamma_{c_{\bf k}}({\bf g})} {\eta}_{c_{\bf k}}({\bf g},{\bf h}).
\end{align}

The vertex basis gauge transformations defined in Eq.~\eqref{eq:vbgt} for quasiparticles generalize naturally to defects.
They are given by
\begin{align}
&\left[ \widetilde{F}^{a_{\bf g} b_{\bf h} c_{\bf k}}_{d_{\bf ghk}} \right]_{(e_{\bf gh}, \alpha,\beta)(f_{\bf hk},\mu,\nu)} \\
&\notag \quad =
\sum_{\alpha', \beta', \mu',\nu'}
[\Gamma^{a_{\bf g} b_{\bf h} }_{e_{\bf gh}}]_{\alpha \alpha'}
[\Gamma^{e_{\bf gh}c_{\bf k}}_{d_{\bf ghk}}]_{\beta \beta'}
\left[ F^{a_{\bf g} b_{\bf h} c_{\bf k}}_{d_{\bf ghk}} \right]_{(e_{\bf gh}, \alpha',\beta')(f_{\bf hk},\mu',\nu')}
[(\Gamma^{b_{\bf h} c_{\bf k}}_{f_{\bf hk} })^{-1} ]_{\mu',\mu}
[(\Gamma^{a_{\bf g} f_{\bf hk}}_{d_{\bf ghk}})^{-1}]_{\nu'\nu}
,
\end{align}
\begin{align}
\left[ \widetilde{R}^{a_{\bf g} b_{\bf h}}_{c_{\bf gh}} \right]_{\mu\nu}
& = \sum_{\mu',\nu'} \left[\Gamma^{b_{\bf h} {}^{\bar{\bf h}}a_{\bf g}}_{c_{\bf gh}} \right]_{\mu \mu'}
\left[R^{a_{\bf g} b_{\bf h}}_{c_{\bf gh}} \right]_{\mu' \nu'} \left[ \left( \Gamma^{a_{\bf g} b_{\bf h}}_{c_{\bf gh}} \right)^{-1} \right]_{\nu',\nu}, \\
\left[\widetilde{U}_{\bf k}(a_{\bf g},b_{\bf h}; c_{\bf gh}) \right]_{\mu \nu} &=
\sum_{\mu',\nu'}
\left[ \Gamma^{{}^{\bf \bar{k}}a_{\bf g} {}^{\bf \bar{k}}b_{\bf h}}_{{}^{\bf \bar{k}}c_{\bf gh}}\right]_{\mu \mu'}
\left[U_{\bf k}(a_{\bf g},b_{\bf h}; c_{\bf gh}) \right]_{\mu' \nu'}
\left[ \left( \Gamma^{a_{\bf g} b_{\bf h}}_{c_{\bf gh}} \right)^{-1} \right]_{\nu',\nu},
\\
\widetilde{\eta}_{c_{\bf k}}({\bf g,h}) & = {\eta}_{c_{\bf k}}({\bf g,h}).
\end{align}
The freedom to add and remove vacuum lines fixes $\gamma_{\I}({\bf h}) = \gamma_a({\bf 0}) =1$ and $\Gamma_{a}^{a\I}=\Gamma_b^{\I b}=\Gamma_{\I}^{\I\I}$.

The $F$-, $R$-, $U$-, and $\eta$-symbols must satisfy the highly constraining pentagon and heptagon equations.
When it is not possible to define a consistent set of this basic data for a given symmetry action and fractionalization class, we say the defectification is obstructed.
This defectification obstruction is captured by an invariant $[\defectO]\in H^4(G,\text{U}(1))$ that is defined in terms of the symmetry action and fractionalization, and which must correspond to the trivial element in order to define a consistent $G$-crossed theory.
When the obstruction vanishes, the defectification classes for a given symmetry action and fractionalization class form an $H^3(G,\text{U}(1))$ torsor.
Beginning from some base theory in the torsor, the distinct theories can be found by gluing in a bosonic SPT characterized by a $3$-cocycle $\alpha({\bf g},{\bf h},{\bf k})\in Z^3(G,\text{U}(1))$.
We present the basic data of bosonic SPT phases, which correspond to $G$-crossed theory whose $\mathcal{C}_{\bf 0}$ contains only the vacuum charge, as an example below.

The full progression of $G$-crossed classification is as follows. (1) Fix the symmetry action $\rho$, specifying how each group element corresponds to a permutation of the topological charges.
(2) Check whether the symmetry action can fractionalize by computing the invariant $[\coho{O}]\in H^3_{[\rho]}(G,\mathcal{A})$.
When $[\coho{O}]$ corresponds to the trivial element, there is no obstruction to fractionalization and the possible symmetry fractionalization classes form an $H^2_{[\rho]}(G,\mathcal{A})$ torsor.
(3) For fixed symmetry action and unobstructed fractionalization class, check whether the theory can consistently be extended to include symmetry defects by computing the invariant $[\defectO]\in H^4(G,\text{U}(1))$.
The theory is unobstructed when $[\defectO]$ corresponds to the trivial element, in which case $\bMTC$ can be extended to a defect theory $\bMTC_G^\times$ characterized by $F$-, $R$-, $U$-, and $\eta$-symbols.
The defectification classes form an $H^3(G,\text{U}(1))$ torsor.
Lastly, we note that the $H^2_{[\rho]}(G,\mathcal{A})$ and $H^3(G,\text{U}(1))$ cohomological classification of $G$-crossed BTCs captures the possible distinct theories that may exist, but potentially overcounts certain ones.
The choice of charge labels is only physical up to relabelings preserving fusion and braiding data.
Thus, the final step of the classification is: (4) check through relabelings and gauge transformations whether any two theories in the torsorial classification can be identified.

\subsection{Example: bosonic symmetry protected topological phases}
\label{sec:bSPT}

The simplest $G$-crosed BTCs describe bosonic SPT phases.
We denote these as $\mathcal{C}_G^\times = \spt_{G}^{[\alpha]}$ with $[\alpha] \in H^3(G,\text{U}(1))$.
Each ${\bf g}$-sector has one simple object, denoted
\begin{align}
\mathcal{C}_{\bf g}= \{\I_{\bf g} \}.
\end{align}
The topological data may be specified for a normalized $3$-cocycle $\alpha \in [\alpha]$ as~\cite{Bark2019}
\begin{align}
F^{\I_{\bf g}\I_{\bf h}\I_{\bf k}} &= \alpha({\bf g}, {\bf h}, {\bf k}), \label{eq:F-alpha}\\
R^{\I_{\bf g}\I_{\bf h}} &= 1 ,\label{eq:R-alpha}\\
U_{\bf k}(\I_{\bf g}, \I_{\bf h}; \I_{\bf gh}) &=
\frac{\alpha({\bf g},{\bf k}, \bar{\bf k}{\bf hk})}{\alpha({\bf g},{\bf h},{\bf k}) \alpha({\bf k},{\bar {\bf k}}{\bf gk},{\bar {\bf k}}{\bf hk})}, \label{eq:U-alpha} \\
\eta_{\I_{\bf k}} ({\bf g},{\bf h})&= \frac{\alpha({\bf g},{\bar {\bf g}}{\bf kg},{\bf h})}{\alpha({\bf g},{\bf h},{\bar{\bf h}\bar{\bf g}}{\bf kgh})\alpha({\bf k},{\bf g},{\bf h})}. \label{eq:eta-alpha}
\end{align}

As noted in the previous section, bosonic SPT phases play a special role in the classification of $G$-crossed theories.
We can obtain a theory $\widehat{\bMTC}_G^\times$ from another theory $\bMTC_G^\times$ with the same symmetry action and fractionalization class by ``gluing'' in $\spt_G^{[\alpha]}$ so that the group labels match those of $\bMTC_G^\times$:
\begin{align}
\label{bSPTgrouplaw}
\widehat{\bMTC}_G^\times &=
\spt_G^{[\alpha]} \underset{G}{\boxtimes} \bMTC_G^\times
\equiv
\spt_G^{[\alpha]} \boxtimes \bMTC_G^\times |_{(\I_{\bf g},a_{\bf g})}.
\end{align}
On the right side, $ |_{(\I_{\bf g},a_{\bf g})}$ indicates that we restrict to the subcategory labeled by objects $(\I_{\bf g},a_{\bf g})$, where ${\bf g} \in G$ and $a_{\bf g} \in \mathcal{C}_{\bf g}$.
The basic data of $\widehat{\bMTC}_{G}^\times$ is given by a product of the data of $\bMTC_G^\times$ and $\spt_G^{[\alpha]}$.
As each theory independently satisfies the pentagon and heptagon equations, the new theory automatically satisfies the $G$-crossed consistency conditions.
Lastly, we note that bosonic SPT phases satisfy a group structure under gluing
\begin{align}\label{eq:SPT-group}
\spt_G^{[\alpha_1]} \underset{G}{\boxtimes} \spt_G^{[\alpha_2]} &= \spt_G^{[\alpha_1 \alpha_2]}
\end{align}
where $[\alpha_1] \cdot [\alpha_2] = [\alpha_1 \alpha_2]$ is the group multiplication law in $H^3(G,\text{U}(1))$.
Equation~\eqref{eq:SPT-group} provides the torsorial action corresponding to the $H^3(G,\text{U}(1))$ part of the classification of $G$-crossed MTCs.

\section{Group Cohomlogy}
\label{app:cohomology}

We briefly review group cohomology; more details can be found in Ref.~\onlinecite{Brown1982}.

Let $G$ be a finite group and $M$ an Abelian group with group action $\rho: G \times M \to M$.
In particular, $\rho$ satisfies
\begin{align}
\rho_{\bf g}(\rho_{\bf h} (a)) &= \rho_{\bf gh}(b)\\
\rho_{\bf g}(a+b) &= \rho_{\bf g}(a) + \rho_{\bf g}(b)
\end{align}
where $a,b \in M$ and we have used additive notation.

Let $\omega({\bf g}_1, {\bf g}_2, \cdots, {\bf g}_n)$ be a function from $n$ group elements valued in $M$.
Such a function is called an `$n$-cochain', and the set of all M-valued $n$-cochains is denoted $C^{n}(G,M)$.
The set of $n$-cochains also forms an Abelian group with multiplication given by,
\begin{align}
[\omega + \omega']({\bf g}_1, \cdots, {\bf g}_n)  =
\omega({\bf g}_1, \cdots, {\bf g}_n)+
\omega'({\bf g}_1, \cdots, {\bf g}_n).
\end{align}
The inverse of $\omega$ is $-\omega$ and the identity element is $\omega({\bf g}_1, \cdots, {\bf g}_n) = 0 \in M$.

The coboundary operator $\text{d}: C^n(G,M) \to C^{n+1}(G,M)$ is defined by
\begin{align}
[\text{d} \omega]({\bf g}_1, \cdots, {\bf g}_{n+1}) = \rho_{{\bf g}_1}[\omega({\bf g}_2, \cdots, {\bf g}_n)] &+ \sum_{j=1}^n (-1)^j \omega({\bf g}_1, \cdots, {\bf g}_{j-1}, {\bf g}_j {\bf g}_{j+1}, {\bf g}_{j+2}, \cdots,{\bf g}_{n+1}) \nonumber \\
&\quad +(-1)^{n+1}\omega({\bf g}_1, \cdots, {\bf g}_{n}).
\end{align}
One can check that $\text{d}\text{d} \omega = 0 \in M$ for any $n$-cochain.

We can now define $n$-cocycles using the coboundary operator $\text{d}$.
An $n$-coycle is given by an $n$-cochain $\omega$ such that $\text{d} \omega = 0$.
The set of $n$-cocycles is given by
\begin{align}
Z^n_\rho(G,M) = \{ \omega \in C^n(G,M) : \text{d} \omega = 0 \}.
\end{align}
If an $n$-cocycle $\omega = \text{d} \mu$, then we say $\omega$ is an $n$-coboundary.
We denote the set of $n$-coboundaries as,
\begin{align}
B^n_\rho(G,M) = \{ \text{d} \mu \in Z^n_\rho(G,M) : \mu \in C^{n-1}(G,M) \}
\end{align}
One can check that $B^n_\rho(G,M) $ is a normal subgroup of $Z^n_\rho(G,M) $.
The $n$-th cohomology group is defined as the quotient of $Z^n_\rho(G,M) $ by $B^n_\rho(G,M) $:
\begin{align}
H^n_\rho (G,M) = \frac{Z^n_\rho(G,M)}{B^n_\rho(G,M)}.
\end{align}

\hbadness=10000	
\bibliography{references}

\end{document}